\documentclass[VANCOUVER,STIX1COL]{WileyNJD-v2}

\articletype{Research Article}%

\received{2020}

\usepackage[font=normalsize]{subcaption}
\usepackage{bm}
\usepackage{ulem}

\newcommand{\spann}{\text{span}}

\usepackage{ulem}

\begin{document}

\title{Non-overlapping block smoothers for the Stokes equations}

\author[1]{Lisa Claus*}

\author[2]{Matthias Bolten**}


\authormark{Claus and Bolten}

\address[1]{*\orgdiv{Computational Research Division}, \orgname{Lawrence Berkeley National Laboratory}, \orgaddress{CA, USA}}

\address[2]{**\orgname{University of Wuppertal}, \orgaddress{\country{Germany}}}


\corres{*\email{lclaus@lbl.gov}\\ **\email{bolten.math@uni-wuppertal.de} }


\abstract[Summary]{Overlapping block smoothers efficiently damp the error contributions from highly oscillatory components within multigrid methods for the Stokes equations but they are computationally expensive. This paper is concentrated on the development and analysis of new block smoothers for the Stokes equations that are discretized on staggered grids. These smoothers are non-overlapping and therefore desirable due to reduced computational costs. Traditional geometric multigrid methods are based on simple pointwise smoothers. However, the efficiency of multigrid methods for solving more difficult problems such as the Stokes equations lead to computationally more expensive smoothers, e.g., overlapping block smoothers. Non-overlapping smoothers are less expensive, but have been considered less efficient in the literature. In this paper, we develop new non-overlapping smoothers, the so-called triad-wise smoothers, and show their efficiency within multigrid methods to solve the Stokes equations. In addition, we compare overlapping and non-overlapping smoothers by measuring their computational costs and
analyzing their behavior by the use of local Fourier analysis.}

\keywords{geometric multigrid methods, Stokes equations, block smoother, local Fourier analysis, finite difference discretization, MAC scheme, staggered grids}


\maketitle


\section{Introduction}\label{sec1}
The Stokes equations are a saddle-point system modeling a slow viscous flow. There is a general interest to find efficient solutions for solving these kind of systems of PDEs, because they describe the physics of many phenomena of scientific and engineering interest \citep{MR2155549}. E.g., the constant movement of the Earth crust caused by the mantle convection is one
application that can be described and simulated by means of the Stokes equations \citep{MR3706420}. Large saddle-point linear systems represent a significant challenge for linear solvers owing to their indefiniteness and often poor spectral properties. Their importance motivates the development of effective solvers.

Multigrid methods have been  applied successfully to large linear systems \citep{Trottenberg_2000, MR1691518, MR4024766}. The methods find an approximate solution through two complementary processes: smoothing and coarse-grid correction. Smoothing filters parts of the error from the approximate solution in a computationally inexpensive way, while coarse-grid correction constructs a lower dimensional problem to remove the error remaining after smoothing. Coarse grid correction consists of the components restriction, interpolation and solution on the coarsest grid. This paper develops efficient multigrid methods to solve the Stokes equations discretized by the MAC scheme. The multigrid components of restriction, interpolation and solution on the coarsest grid can be immediately generalized to large linear systems of saddle point type.

However, the design of efficient smoothers for solving systems of PDEs requires special attention. We distinguish between decoupled smoothing and collective smoothing. A zero diagonal block that appears in the discrete system hampers a basic numerical treatment of
the problem, which would be to relax the discrete equations directly in a decoupled way. A variety of multigrid methods  for the Stokes equations based on different smoothing approaches have been developed, used and analyzed \citep{MR3580776, MR3323212, MR3805850, MR2001083}. We focus on collective relaxation, so-called block smoothers. A comparison of a well-known overlapping block smoother, the Vanka smoother, with a non-overlapping smoother, the triad-wise smoother, is given. While the Vanka smoother is well established to solve the Stokes equations it is computationally expensive \citep{VANKA1986138, clausphd}. Non-overlapping block smoothers are less expensive but were considered less efficient in terms of their smoothing and convergence behavior \citep{VANKA1986138}. In this paper, we develop a new non-overlapping block smoother based on the results of the comparison with both smoothers. For periodic boundary conditions we show the diminished efficiency of the triad-wise smoother relative to the Vanka smoother. Solving the Stokes equations with Dirichlet boundary conditions leads to an enormous decrease of the efficiency for the triad-wise smoother. We highlight the issues of the non-overlapping smoother by means of analysis and numerical examples and demonstrate the causes of the high computational costs of the overlapping smoother. Based on these results a new non-overlapping smoother is constructed for Dirichlet boundary conditions that overcomes the demonstrated issues. To the best of our knowledge, the proposed modified smoother represents the first-ever non-overlapping block smoother that reaches efficiency close to the overlapping block smoother.

Local Fourier analysis (LFA) is utilized as the main analysis tool for quantitative convergence estimates and to optimize geometric multigrid components such as smoothers or intergrid transfer operators. The idea was introduced by Brandt \citep{BRANDT1977} and later extended and refined in \citep{BRANDT1994}. Practical guidelines on the use of LFA can be found in \citep{MR2108045}. LFA is commonly applied to relaxation schemes using (symmetric) Gauss-Seidel approaches \citep{Niestegge_1990}, and has been extented to more sophisticated relaxation schemes \citep{MR3795547, MR3922504, MR3488076}.  In this study, LFA is used to analyze the two-grid method based on the overlapping and non-overlapping smoothers. In addition to the LFA results, we provide numerical results that guide the development of a new algorithm leading to more efficient multigrid methods.

The outline of the paper is as follows. In Section 2, we introduce the marker-and-cell (MAC) finite-difference discretization of the two-dimensional Stokes equations and the basic idea of the two block relaxation methods used in this paper. In Section 3, we present some basic definitions and notations to introduce LFA for the Stokes equations. In addition, analysis results are given at the end of the section. In Section 4, numerical examples extend the analysis and highlight the differences and issues arising due to different boundary conditions. In Section 5, we develop a new non-overlapping block smoother, an update of the so-called triad-wise smoothers leading to an enormous increase of efficiency. We give some details about parallelization aspects in Section 6. Conclusions are drawn in Section 7.

\section{Smoother for Stokes equations}\label{sec2}

\subsection{Stokes equations}
\begin{figure}[h]
\centering
{\includegraphics[width=4 cm]{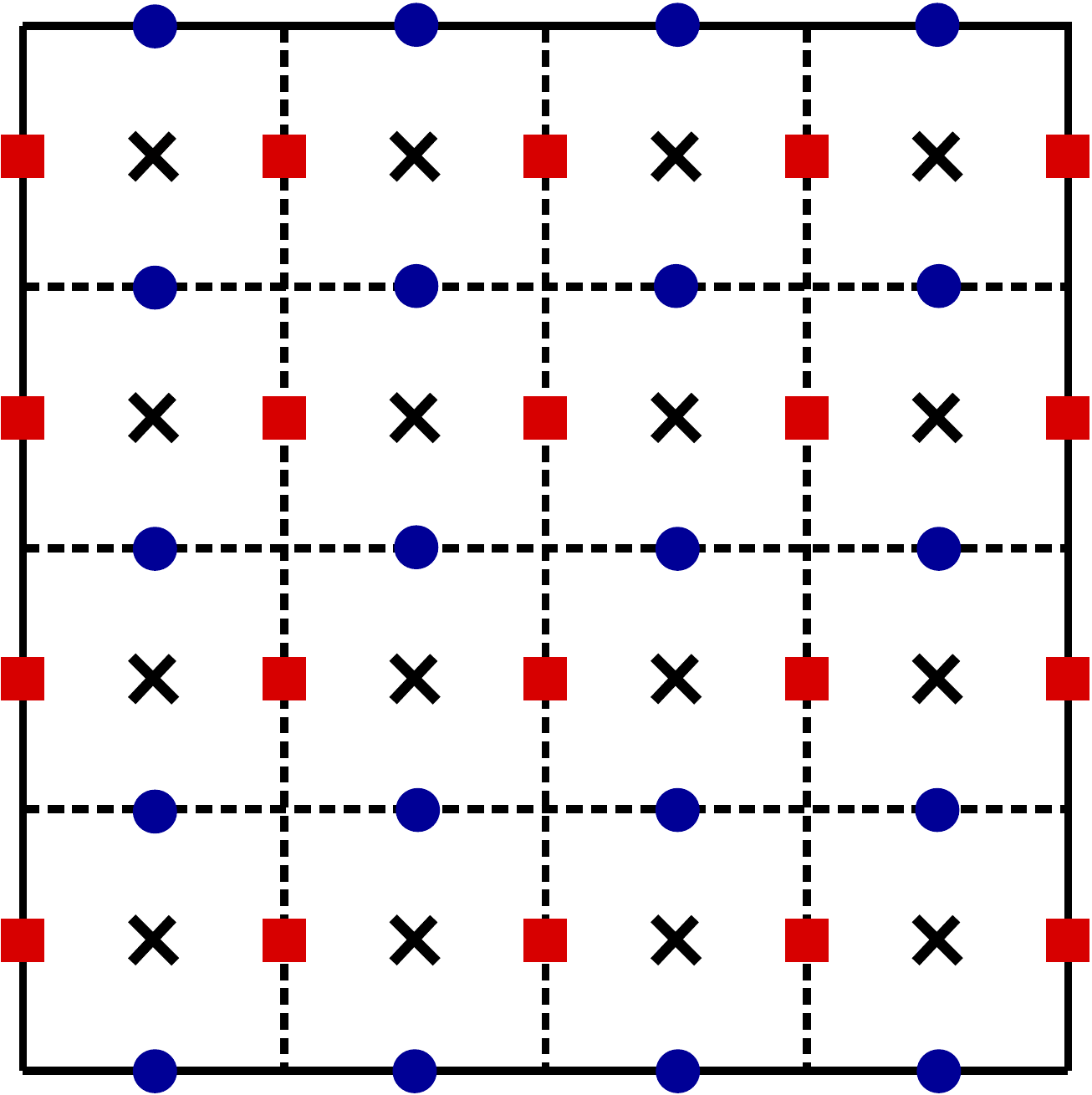}}
\caption[]{Location of unknowns $(u,v,p)$.}
\label{fig:1}
\end{figure}

We introduce the Stokes equations in two dimensions here. Let $\Omega$ be a domain in $\mathbb{R}^2$ with boundary $\partial \Omega$, and $\mathbf{f}$ be a given vector function that describes an external force. Solving the two-dimensional Stokes equations means finding a fluid velocity $\mathbf{u} = (u, v)$ and pressure $p$ such that the following system is fulfilled,
\begin{align}
\Delta u + \partial_{\mathrm{x}_1} p &= \mathrm{f}_1,\label{eq:1}\\
\Delta v + \partial_{\mathrm{x}_2} p &= \mathrm{f}_2,\label{eq:2}\\
\partial_{\mathrm{x}_1} u + \partial_{\mathrm{x}_2} v &= 0,\label{eq:3}
\end{align}
where $\partial_{\mathrm{x}_i}:=\frac{\partial}{\partial {\mathrm{x}_i}}$.
The discretization of the Stokes equations naturally leads to linear systems of the form
\begin{align}
\label{eq:4}
L_h=
\begin{pmatrix}
A_1 & & B_1\\
& A_1 & B_2\\
-B_1^T & -B_2^T & 
\end{pmatrix} \begin{pmatrix}
u_h \\ v_h \\ p_h
\end{pmatrix} = \begin{pmatrix}
\mathrm{f}_1 \\ \mathrm{f}_2 \\ 0
\end{pmatrix},
\end{align}
so-called saddle point systems. Here $A_1$ denotes the negative discrete Laplace operator, $A_1=-\Delta_h$, and $B_1^T$ and $B_2^T$ are the discrete divergence operators, $B_1^T=(\partial_{\mathrm{x}_1})_{h/2}$ and $B_2^T=(\partial_{\mathrm{x}_2})_{h/2}$. 
In this paper, we use a staggered distribution of unknowns, the so-called Marker-and-Cell (MAC) scheme. The idea of the MAC scheme is to place the unknowns ($u, v, p$) in different locations. More specifically, the discrete pressure unknowns $p$ are located in the center of each cell and the discrete values of the velocity unknowns are located on the midpoint of vertical edges (${\mathrm{x}_1}$-component $u$) and on the midpoint of horizontal edges (${\mathrm{x}_2}$-component $v$), as seen in Figure \ref{fig:1}. We consider uniform meshes with: $h_x = h_y = h$ in this paper.

The discrete analogue of the ${\mathrm{x}_1}$-coordinate momentum equation \eqref{eq:1} is defined at vertical edges, the discrete ${\mathrm{x}_2}$-coordinate momentum equation \eqref{eq:2} at horizontal edges, and the continuity equation \eqref{eq:3} at cell centers using central difference schemes. We introduce the following indexing system with regard to the underlying grid: $i$ is the column index and $j$ is the row index, ranging from $1:n$ or $1:n+1$, where $n$ is the number of cells in one direction. The discrete pressure unknowns are $p_{i,j}$, where $i=1,\dots,n$, $j=1,\dots,n$, and the velocity unknowns are $u_{i,j}$ with $i=1,\dots,n+1$, $j=1,\dots,n$ and $v_{i,j}$ with $i=1,\dots,n$, $j=1,\dots,n+1$. Using this indexing system, the MAC scheme can be written as
\begin{align}
\frac{4u_{i,j} - u_{i-1,j} - u_{i+1,j} - u_{i,j-1} - u_{i,j+1}}{h^2} + \frac{p_{i,j}-p_{i-1,j}}{h}&=f^{i,j}_1, \label{eq:5}\\
\frac{4v_{i,j} - v_{i-1,j} - v_{i+1,j} - v_{i,j-1} - v_{i,j+1}}{h^2} + \frac{p_{i,j}-p_{i,j-1}}{h}&=f^{i,j}_2, \label{eq:6}\\
\frac{u_{i+1,j}-u_{i,j}}{h} + \frac{v_{i,j+1}-v_{i,j}}{h} &= 0. \label{eq:7}
\end{align}

Our interest in this paper is to develop efficient multigrid methods to
solve the Stokes equations discretized by the MAC scheme. We discuss how the multigrid components of smoothing, restriction, interpolation and solution on the coarsest grid are generalized to saddle point systems. A key point of an efficient multigrid method is the right choice of the relaxation procedure. The design of efficient smoothers for solving systems of PDEs often requires special attention. The relaxation method should smooth the error for all unknowns in the equations, which is not an easy task. Therefore, we focus on constructing appropriate relaxation schemes to smooth the error of the unknown functions. The other multigrid components can immediately be extended to systems of PDEs.
A zero diagonal block that appears in the discrete system hampers a basic numerical treatment of the problem, which would be to relax the discrete equations directly in a decoupled way. More specifically, decoupled iterative solution methods use blocks on the main diagonal in an equation-wise fashion. Unfortunately, this equation-wise relaxation is not possible due to the zero diagonal block. Moreover, a first obvious choice for these so-called strong off-diagonal operators in the differential system is collective smoothing. However, for staggered grids the unknowns
are not defined at the same locations and therefore we define standard collective relaxation as described in the following section.

\subsection{Vanka smoother}

\begin{figure}[h]
\centering
\subcaptionbox[]{All five unknowns of one block
    \label{subfig:2a}} 
[0.45\textwidth] 
{\includegraphics[width=3 cm]{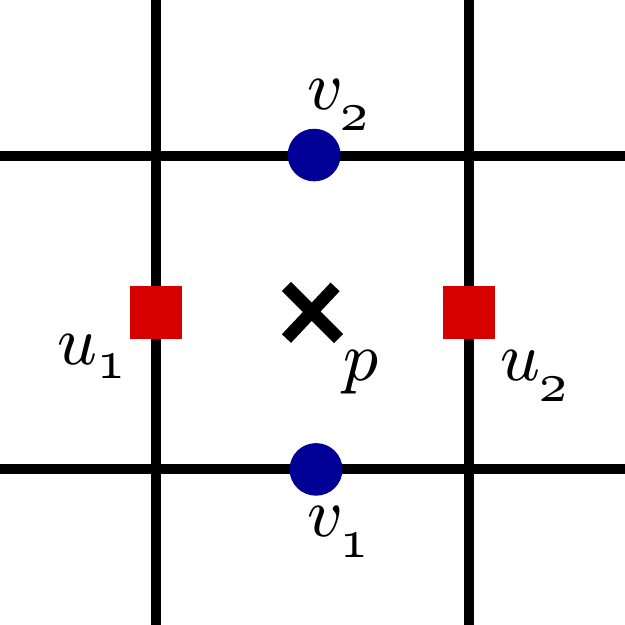}}
\subcaptionbox[]{Two overlapping cells
    \label{subfig:2b}} 
[0.45\textwidth] 
{\includegraphics[width=4.3 cm]{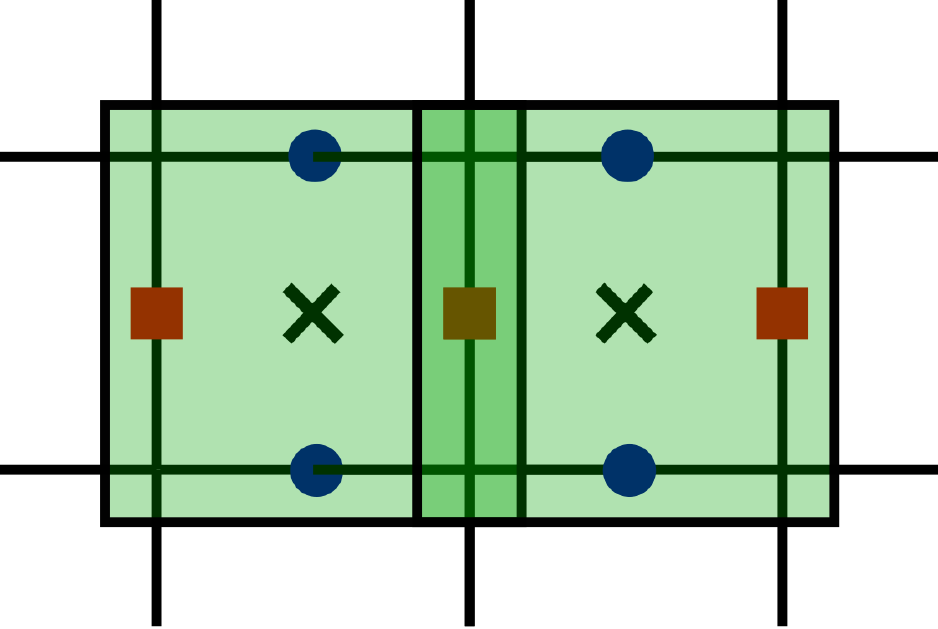}}
\caption[]{The Vanka smoother: five unknowns updated collectively.}
\label{fig:2}
\end{figure}

The basic idea of the Vanka smoother is to solve the discrete Stokes equations locally ``cell by cell''. This means all five unknowns of one cell are updated collectively, involving the respective four momentum equations at the cell boundaries and one continuity equation in the center of the cell, c.f.\ Figure \ref{fig:2} and Equations \eqref{eq:5}-\eqref{eq:7}. This results in an overlapping updating process of $5 \times 5$ systems of equations. This method was first introduced in \citep{VANKA1986138}. 

Using this scheme, each velocity component is updated twice and the pressure once per relaxation. The process of relaxation appears in a Gauss-Seidel manner. That means, after smoothing all unknowns of one cell collectively, we use the updated values for the relaxation of the next cell. This relaxation method is commonly referred to as \textit{overlapping Schwarz smoother} or \textit{symmetric coupled Gauss-Seidel} (SCGS). Especially in the context of the Stokes equations, we use the term \textit{Vanka smoother}. It turns out to have robust smoothing properties \citep{Trottenberg_2000}. Due to the overlap the iteration matrix $M_h$ of the relaxation scheme cannot be easily represented by a banded matrix or equivalently as a stencil. This is possible if a suitable coloring scheme is chosen. More details can be found in Section \ref{sec6}.

\subsection{Triad-wise smoother}

For the Stokes system, the Vanka smoother is a very  efficient smoother \citep{VANKA1986138}. However, due to the overlap, computational costs are high and parallel implementation aspects are not satisfying due to the required coloring scheme, c.f. Section \ref{sec6}. We exploit the potential of the Vanka relaxation method and focus on the development of a block smoother that is less expensive and still efficient.

\begin{figure}[h]
\centering
\subcaptionbox[]{All three unknowns of one triad block
    \label{subfig:3a}} 
[0.45\textwidth] 
{\includegraphics[width=3 cm]{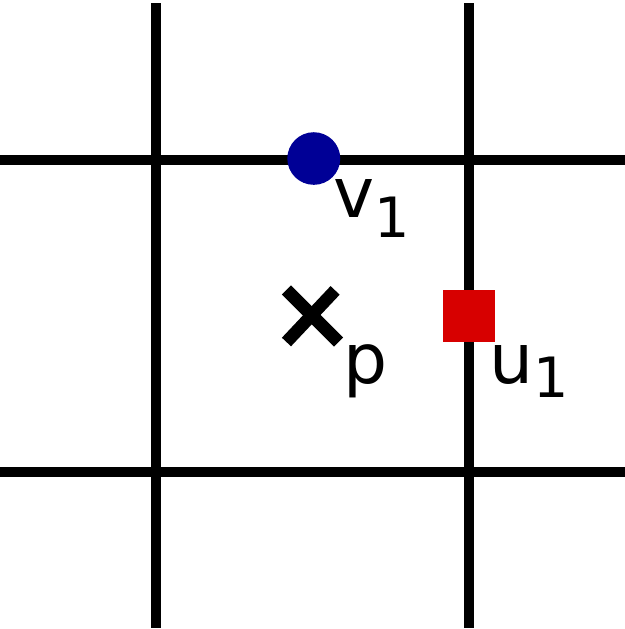}}
\subcaptionbox[]{Two non-overlapping blocks
    \label{subfig:3b}} 
[0.45\textwidth] 
{\includegraphics[width=4.3 cm]{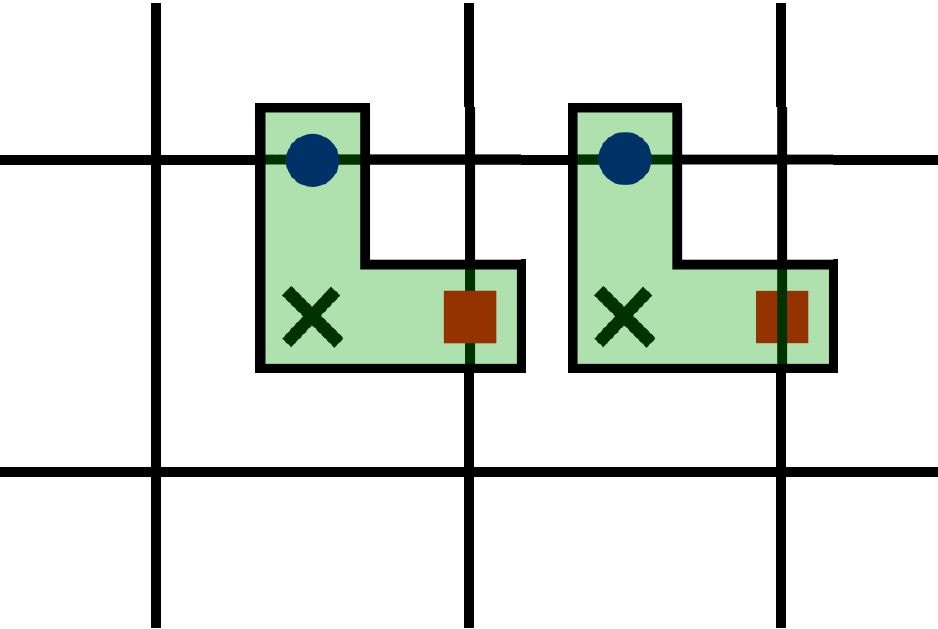}}
\caption[]{Unknowns updated collectively by triad-wise relaxation.}
\label{fig:3}
\end{figure}

The basic idea of the triad-wise smoother is to update three unknowns collectively, involving two momentum equations at cell boundaries and the continuity equation in the center of the cell. This results in a non-overlapping updating process of $3 \times 3$ systems of equations, as seen in Figure \ref{fig:3}. Using this scheme, one velocity component in each direction and one pressure component are updated simultaneously. The first version of the triad-wise smoother was introduced in \citep{VANKA1986138} and is referred to as \textit{unsymmetric coupled Gauss-Seidel} (UCGS) \citep{VANKA1986138}. The potential of a different updating process as for example a Jacobi-manner or a different ordering instead of the Gauss-Seidel updating process will be examined in this paper. We stick to the term triad-wise smoother whenever we don't specify the manner of the updating process.
The interation matrix of the relaxation scheme, $M_h$ can be represented by
\begin{align}
\label{eq:8}
\begin{split}
M_h:= \frac{1}{\omega} D_h-L_h = \begin{bmatrix}
-\frac{1}{h^2} & & -\frac{1}{h}\\
&-\frac{1}{h^2} & \\
& &
\end{bmatrix} &\begin{bmatrix}
\textcolor{white}{+}\frac{4}{h^2} & & \textcolor{white}{+}\frac{1}{h}\\
&\textcolor{white}{+}\frac{4}{h^2} & \textcolor{white}{+}\frac{1}{h} \\
\textcolor{white}{+}\frac{1}{h}& \textcolor{white}{+}\frac{1}{h} &
\end{bmatrix},\\
&\begin{bmatrix}
-\frac{1}{h^2} & & \\
&-\frac{1}{h^2} & -\frac{1}{h}\\
& &
\end{bmatrix}
\end{split}
\end{align} for a Gauss-Seidel-type updating process and 
\begin{align}
\label{eq:9}
M_h:= \frac{1}{\omega} D_h = \frac{1}{\omega} \begin{bmatrix}
\textcolor{white}{+}\frac{4}{h^2} & & \textcolor{white}{+}\frac{1}{h}\\
&\textcolor{white}{+}\frac{4}{h^2} & \textcolor{white}{+}\frac{1}{h} \\
\textcolor{white}{+}\frac{1}{h}& \textcolor{white}{+}\frac{1}{h} &
\end{bmatrix},
\end{align}
for Jacobi-type updating steps.
In the light of higher computational cost, the triad-wise smoother seems to be a good alternative to the Vanka smoother, since it has no overlap and therefore reduces computational costs. In addition, this relaxation method has better parallelization properties. These characteristics are examined in detail in the next section. However, the unsymmetric coupled (triad-wise) smoother was observed to have poor smoothing and often led to divergence for the Stokes equations \citep{VANKA1986138}.

In the following we analyze and describe the conditions that lead to poor smoothing of the triad-wise smoother. In addition, we investigate different versions of the triad-wise smoother and modify the relaxation method to generate an efficient smoother with good smoothing properties.

\section{Local Fourier analysis for Stokes equations}\label{sec3}

We employ local Fourier analysis (LFA) to analyze the convergence properties of multigrid methods for the Stokes equations and thereby compare the triad-wise and Vanka smoother. 

\subsection{Definitions and notations}\label{sec31}

In order to describe LFA for staggered grids, we first introduce some terminology. More details can be found, for example, in \cite{Trottenberg_2000}. We consider two-dimensional infinite uniform grids $G_h = G^1_h \cup G^2_h \cup G^3_h$ where
\begin{align}
\label{eq:10}
G_h^j = \{ \bm{x}^j := \bm{kh} + \delta^j, \bm{k} \in \mathbb{Z}^2 \}, \text{with } \delta^j = \begin{cases} (0, h/2) & \text{if } j=1,\\ 
(h/2, 0) & \text{if } j=2,\\
(h/2, h/2) & \text{if } j=3. \end{cases}
\end{align}
Here, the grids that correspond to the velocity unknowns $u$ and $v$,
are denoted by $G_h^1$ and $G_h^2$. Moreover, the grid that corresponds to pressure unknown $p$ is denoted by $G_h^3$. The coarse grid, $G_H$, is constructed by doubling the mesh size, i.e. we use standard coarsening. The functions $\varphi(\bm{\theta, x})$ are called Fourier modes or Fourier functions with frequency $\bm{\theta}$ and form the basis of the Fourier space $\mathcal{F}_h$. 
Under general assumptions, each multigrid operator can be replaced by an operator with constant coefficients. Therefore, Fourier modes serve as eigenfunctions of multigrid operators and can be used to analyze the effect of each multigrid component on different modes.
For staggered grids, the basis of the Fourier space $\mathcal{F}_h$ is given by the $\spann \{ \varphi^1(\bm{\theta, x}), \varphi^2 (\bm{\theta, x}), \varphi^3 (\bm{\theta, x})\}$ for $\bm{\theta} \in \left[-\frac{\pi}{2},\frac{3\pi}{2}\right)^2$, with
\begin{align*}
\varphi^1(\bm{\theta, x})&= \begin{pmatrix}
e^{\imath \bm{\theta} \bm{x}^1/\bm{h}} & 0 & 0 
\end{pmatrix}^T, \varphi^2(\bm{\theta, x})= \begin{pmatrix}
0 & e^{\imath \bm{\theta} \bm{x}^2/\bm{h}} & 0 
\end{pmatrix}^T,\\
\varphi^3(\bm{\theta, x})&= \begin{pmatrix}
0 & 0 & e^{\imath \bm{\theta} \bm{x}^3/\bm{h}} 
\end{pmatrix}^T.
\end{align*}
We use the stencil notation and consider the operator $L_h$ of the Stokes system,
\begin{align*}
L_h=\begin{pmatrix}
-\Delta_h & & (\partial_{\mathrm{x}_1})_{h/2}\\
& -\Delta_h & (\partial_{\mathrm{x}_2})_{h/2}\\
-(\partial_{\mathrm{x}_1})_{h/2} & - (\partial_{\mathrm{x}_2})_{h/2}& \\
\end{pmatrix},
\end{align*}
with stencils 
\begin{align*}
-\Delta_h=\frac{1}{h^2}\begin{bmatrix}
& -1 & \\
-1 & 4 & -1\\
& -1 & \\
\end{bmatrix}, \quad
(\partial_{\mathrm{x}_1})_{h/2} = \frac{1}{h} \begin{bmatrix}
-1 & 0 & 1
\end{bmatrix}, \quad 
(\partial_{\mathrm{x}_2})_{h/2} = \frac{1}{h} \begin{bmatrix}
1 \\ 0 \\ -1
\end{bmatrix}.
\end{align*}
\begin{definition}
We call $\tilde{L}_h$ the symbol of $L_h$.
\end{definition}
Note that for all Fourier modes $\varphi(\bm{\theta, x})$,
\begin{align*}
 L_h \varphi(\bm{\theta, x}) = \tilde{L}_h \varphi(\bm{\theta, x}).
\end{align*}
We assume that the systemmatrix $L_h$ is a block matrix, where each block is an infinite-grid Toeplitz matrix. Then, each block in $L_h$ may be diagonalized with LFA.  Each entry in the symbol $\tilde{L}_h$ is computed as the (scalar) symbol of the corresponding block of $L_h$. Since $L_h$ is a $3 \times 3$ block operator, its symbol is a $3 \times 3$ matrix,
\begin{align}
\label{eq:11}
\tilde{L}_h(\bm{\theta})=\begin{pmatrix}
\frac{1}{h^2}(4-2\cos\theta_1-2\cos\theta_2) & 0 & 2\imath h \sin \frac{\theta_1}{2}\\
0 & \frac{1}{h^2}(4-2\cos\theta_1-2\cos\theta_2) & 2 \imath h \sin \frac{\theta_2}{2}\\
-2\imath h \sin \frac{\theta_1}{2} & -2\imath h \sin \frac{\theta_2}{2}& 0 \\
\end{pmatrix}.
\end{align}
High and low frequencies for standard coarsening are given by
\begin{align*}
\bm{\theta} \in T^{\text{low}}:=\left[-\frac{\pi}{2}, \frac{\pi}{2}\right)^2, \quad \bm{\theta} \in T^{\text{high}}:=\left[-\frac{\pi}{2},\frac{3\pi}{2}\right)^2\Big{\backslash} \left[-\frac{\pi}{2}, \frac{\pi}{2}\right)^2.
\end{align*} 
\begin{remark} 
The influence of boundary conditions is not taken into account when applying LFA, since the modes $\varphi_h(\bm{\theta,x})$ are defined on the infinite grid $G_h$. Experience with LFA show that it often serves as an accurate prediction tool for problems with periodic boundary conditions, but degradation in performance may be seen with Dirichlet boundary conditions \citep{Niestegge_1990}.
\end{remark}
We start with a two-grid analysis and compare results for both smoothers.

\subsection{Two-grid LFA}

We analyze the two-grid method by starting with the analysis aspects of the smoothing procedures. This is followed by an introduction of coarse-grid correction symbols.\\
Due to the overlapping cells of the Vanka smoother, certain degrees of freedom are updated multiple times over the course of a single sweep of relaxation. Therefore, classical LFA techniques fail. LFA is based on the assumption that the iteration matrix $M_h$ can be written as a Toeplitz operator. It is not apparent that the same is true for overlapping smoothers.\\
An extension of the classical LFA analysis for the Vanka smoothers was first introduced by Sivaloganathan in \citep{SIVALOGANATHAN1991246}. MacLachlan and Oosterlee generalized the LFA techniques to analyze overlapping smoothers for other PDEs and discretizations \citep{MacLachlanOosterlee2011}. Their paper provides theoretical results in order to use the Fourier ansatz. A necessary assumption is that the error-propagation operator for coupled (overlapping) relaxation is block-diagonalized by the Fourier matrix, regardless of the distribution of degrees of freedom within the block. The key step in proving this is to show the following inductive step: if the errors before relaxation on degrees of freedom within a block $V_{i,j}$ satisfy a generalized LFA ansatz, then the errors after relaxation on $V_{i,j}$ satisfy the same ansatz advanced by one cell. MacLachlan and Oosterlee show that the error-propagation matrix for any coupled relaxation is an infinite-grid block-multilevel-Toeplitz matrix \citep{MacLachlanOosterlee2011}. This means that we can attempt to analyze these techniques using classical multigrid smoothing and two-grid Fourier analysis tools to measure the effectiveness of the resulting multigrid cycles.\\
Practical guidelines to obtain the smoothing symbol $\tilde{\mathcal{S}}_h(\bm{\theta})$ of the Vanka smoother to analyze this relaxation method are provided in \citep{clausphd}. For more details and a generalization to other overlapping smoothers, we refer to \citep{MacLachlanOosterlee2011}.\\
Due to the non-overlapping blocks, the LFA analysis for the triad-wise smoother is based on a simple expansion of the assumptions of LFA for scalar PDEs as described in Section \ref{sec31}. The error-propagation operator for triad-type smoothers, represented by an operator $M_h$ is
\begin{align*}
\mathcal{S}_h(\omega)=I- M_h^{-1}L_h.
\end{align*} 
The corresponding smoothing symbol $\tilde{\mathcal{S}}_h(\omega,\bm{\theta})$ is given by
\begin{align*}
\tilde{\mathcal{S}}_h(\omega, \bm{\theta})=I- \tilde{M}_h^{-1}(\bm{\theta})\tilde{L}_h(\bm{\theta}),
\end{align*}
with $I$ the identity matrix and $\tilde{L}_h(\bm{\theta})$ as given in \eqref{eq:11}. The symbol of the block iteration matrix $\tilde{M}_h$ can be given by
\begin{align*}
\tilde{M}_h(\bm{\theta})= \begin{pmatrix}
\frac{1}{h^2}\left(\frac{4}{\omega}-e^{-\imath \theta_1}-e^{-\imath \theta_2}\right)& 0 & \frac{1}{h} \left(\frac{1}{\omega}e^{\frac{1}{2}\imath \theta_1}-e^{-\frac{1}{2}\imath \theta_1}\right)\\
0 & \frac{1}{h^2}\left(\frac{4}{\omega}-e^{-\imath \theta_1}-e^{-\imath \theta_2}\right)& \frac{1}{h} \left(\frac{1}{\omega}e^{\frac{1}{2}\imath \theta_2}-e^{-\frac{1}{2}\imath \theta_2}\right)\\
\frac{1}{\omega h}e^{-\frac{1}{2}\imath \theta_1}& \frac{1}{\omega h}e^{-\frac{1}{2}\imath \theta_2} & 0
\end{pmatrix},
\end{align*}
for a weighted Gauss-Seidel updating process and by
\begin{align*}
\tilde{M}_h(\bm{\theta})= \frac{1}{\omega}\begin{pmatrix}
\frac{4}{h^2}& 0 & \frac{1}{h} e^{\frac{1}{2}\imath \theta_1}\\
0 & \frac{4}{h^2}& \frac{1}{h} e^{\frac{1}{2}\imath \theta_2}\\
\frac{1}{h}e^{-\frac{1}{2}\imath \theta_1}& \frac{1}{h}e^{-\frac{1}{2}\imath \theta_2} & 0
\end{pmatrix},
\end{align*}
for a weighted Jacobi updating process.\\

For systems of PDEs that are discretized on staggered grids, the analysis of the coarse grid correction is challenging due to the different location of unknowns. In that case, LFA of transfer operators cannot be done as in the scalar case. MacLachlan and Oosterlee \citep{MacLachlanOosterlee2011} discussed Fourier representations of grid-transfer operators for general staggered meshes in the context of systems of PDEs. Here, we have three types of grid points on the fine and coarse grid. The restriction operator can be decomposed based on the partitioning of dofs associated with the location of the unknowns on the grid. Osterlee and MacLachlan introduced, explained and proved the need to modify LFA due to this staggering. We emphasize this statement in Remark \ref{rem2}. More details can be found in \citep{MacLachlanOosterlee2011, Yunhui2018TWOLEVELFA, MR3922504}. 
\begin{remark} \label{rem2} 
When calculating the symbol of a restriction operator that mixes different types of dofs, we must split it into the different types of dofs that it restricts from and to. If the restriction
operator is defined on a staggered grid, we have $G_H= G^1_H \cup G^2_H \cup G^3_H$ (c.f.\ \eqref{eq:10}) and the symbols depend on the location $\bm{x}$ of the unknowns.
\end{remark}
In the next paragraph, we show how to determine the symbol of a restriction operator defined on a staggered mesh. A similar structure for interpolation is established subsequently. Let $\varphi_h(\bm{\theta^{\xi}, x}) = e^{\imath \bm{\theta^{\xi}} \bm{ x}/\mathbf{h}}$. We have the following equality
\begin{align*}
\varphi_h (\bm{\theta^{\xi}, x}) = e^{\imath \bm{\xi} \pi \bm{x}/\mathbf{h}} \varphi_H(2\bm{\theta^{(0,0)}, x}), \quad \text{for all } \bm{x} \in G_H.
\end{align*}
Note that this equality is different to the classical theory and was developed by MacLachlan and Oosterlee in \citep{MacLachlanOosterlee2011}. The following relations are based on this update of the classical theory. Due to the definition of $G_H= G^1_H \cup G^2_H \cup G^3_H$, we have
\begin{align*}
\bm{x}=
\begin{cases}
\begin{pmatrix} 2jh, & 2(l+\frac{1}{2})h \end{pmatrix}^T, & \text{for } j,l \in \mathbb{Z}, \text{ if } \bm{x} \in G_H^1,\\
\begin{pmatrix} 2(j+\frac{1}{2})h, & 2lh \end{pmatrix}^T, & \text{for } j,l \in \mathbb{Z}, \text{ if } \bm{x} \in G_H^2,\\
\begin{pmatrix} 2(j+\frac{1}{2})h, & 2(l+\frac{1}{2})h \end{pmatrix}^T, & \text{for } j,l \in \mathbb{Z}, \text{ if } \bm{x} \in G_H^3.
\end{cases}
\end{align*}
Note that the following relations hold,
\begin{align*}
e^{\imath \bm{\xi} \pi \bm{x}/\mathbf{h}}=
\begin{cases}
e^{\imath \xi_1 \pi 2j} e^{\imath \xi_2 \pi 2l} e^{\imath \xi_2 \pi} = e^{\imath \pi \xi_2} = (-1)^{\xi_2}, & \text{for } j,l \in \mathbb{Z}, \text{ if } \bm{x} \in G_H^1,\\
e^{\imath \xi_1 \pi 2j} e^{\imath \xi_1 \pi} e^{\imath \xi_2 \pi 2l}  = e^{\imath \pi \xi_1} = (-1)^{\xi_1}, & \text{for } j,l \in \mathbb{Z}, \text{ if } \bm{x} \in G_H^2,\\
e^{\imath \pi \xi_1} e^{\imath \pi \xi_2} = (-1)^{\xi_1} (-1)^{\xi_2}, & \text{for } j,l \in \mathbb{Z}, \text{ if } \bm{x} \in G_H^3.
\end{cases}
\end{align*}
Then, the Fourier representation of $R_h^H$ is given by the $(3 \times 12)$ matrix
\begin{align*}
\hat{R}_h^H(\bm{\theta})&=\begin{pmatrix}
\tilde{R}_h^H(\bm{\theta}^{(0,0)}) & \tilde{R}_h^H(\bm{\theta}^{(1,1)}) & \tilde{R}_h^H(\bm{\theta}^{(1,0)}) & \tilde{R}_h^H(\bm{\theta}^{(0,1)})
\end{pmatrix},
\end{align*}
where $\tilde{R}_h^H(\bm{\theta^\xi})$ is given in Definition \ref{def2}. Note that the matrices $\tilde{R}_h^H(\bm{\theta^\xi})$ are diagonal matrices since the restriction operators of different type of unknowns are not coupled,
\begin{align*}
\tilde{R}_h^H(\bm{\theta^\xi}) = \begin{pmatrix}
\tilde{R}_h^H(\bm{\theta^\xi})\big\lvert_u & & \\
& \tilde{R}_h^H(\bm{\theta^\xi})\big\lvert_v & \\
& & \tilde{R}_h^H(\bm{\theta^\xi})\big\lvert_p 
\end{pmatrix}.
\end{align*}
\begin{definition}
\label{def2}
We call $\tilde{R}_h^H(\bm{\theta^\xi}) := \sum\limits_{\bm{\kappa} \in \mathbf{V}} \mathrm{r}_{\bm{\kappa}} e^{\imath \bm{\theta^\xi} \cdot \bm{\kappa}} e^{\imath \bm{\xi} \pi \bm{x}/\mathbf{h}}$ the symbol of $R_h^H$.
\end{definition} 
\begin{definition} \label{def3} 
The restriction operators for the Stokes equations on staggered grids are given by
\begin{align*}
\renewcommand\arraystretch{1.3}
R_h^H\big\lvert_u=\begin{bmatrix}  
\frac{1}{8} & \frac{1}{4} & \frac{1}{8} \\   
& \ast & \\   
\frac{1}{8} & \frac{1}{4} & \frac{1}{8} \end{bmatrix}, \quad R_h^H\big\lvert_v=\begin{bmatrix}  
\frac{1}{8} & & \frac{1}{8} \\ 
\frac{1}{4} & \ast & \frac{1}{4} \\   
\frac{1}{8} & & \frac{1}{8} \end{bmatrix},  \quad
R_h^H\big\lvert_p=\begin{bmatrix}  
\frac{1}{4} & & \frac{1}{4} \\   
& \ast & \\
\frac{1}{4} & &\frac{1}{4} \end{bmatrix},
\end{align*}
where $\ast$ denotes the resulting coarse-grid point.
\end{definition}

Based on Definition \ref{def3}, we are able to compute the entries of $\hat{R}_h^H(\bm{\theta})$. We show the computation of the symbol for the restriction operator $R_h^H\big\lvert_u$, i.e., the first row of the matrix $\hat{R}_h^H(\bm{\theta})$. 
{\footnotesize{\begin{align*}
\tilde{R}_h^H(\bm{\theta}^{(0,0)})\big\lvert_u &= \frac{1}{8} e^{-\imath \theta_1} e^{-\frac{1}{2}\imath \theta_2} + \frac{1}{4} e^{-\frac{1}{2}\imath \theta_2} + \frac{1}{8} e^{\imath \theta_1} e^{-\frac{1}{2}\imath \theta_2} + \frac{1}{8} e^{-\imath \theta_1} e^{\frac{1}{2}\imath \theta_2} + \frac{1}{4} e^{\frac{1}{2}\imath \theta_2} + \frac{1}{8} e^{\imath \theta_1} e^{\frac{1}{2} \imath \theta_2},\\
\tilde{R}_h^H(\bm{\theta}^{(1,1)})\big\lvert_u &= -\frac{1}{8} e^{-\imath \theta_1} e^{-\frac{1}{2}\imath \theta_2} - \frac{1}{4} e^{-\frac{1}{2}\imath \theta_2} - \frac{1}{8} e^{\imath \theta_1} e^{-\frac{1}{2}\imath \theta_2} - \frac{1}{8} e^{-\imath \theta_1} e^{\frac{1}{2}\imath \theta_2} - \frac{1}{4} e^{\frac{1}{2}\imath \theta_2} - \frac{1}{8} e^{\imath \theta_1} e^{\frac{1}{2} \imath \theta_2},\\
\tilde{R}_h^H(\bm{\theta}^{(1,0)})\big\lvert_u &= \frac{1}{8} e^{-\imath \theta_1} e^{-\frac{1}{2}\imath \theta_2} + \frac{1}{4} e^{-\frac{1}{2}\imath \theta_2} + \frac{1}{8} e^{\imath \theta_1} e^{-\frac{1}{2}\imath \theta_2} + \frac{1}{8} e^{-\imath \theta_1} e^{\frac{1}{2}\imath \theta_2} + \frac{1}{4} e^{\frac{1}{2}\imath \theta_2} + \frac{1}{8} e^{\imath \theta_1} e^{\frac{1}{2} \imath \theta_2},\\
\tilde{R}_h^H(\bm{\theta}^{(0,1)})\big\lvert_u &= -\frac{1}{8} e^{-\imath \theta_1} e^{-\frac{1}{2}\imath \theta_2} - \frac{1}{4} e^{-\frac{1}{2}\imath \theta_2} - \frac{1}{8} e^{\imath \theta_1} e^{-\frac{1}{2}\imath \theta_2} - \frac{1}{8} e^{-\imath \theta_1} e^{\frac{1}{2}\imath \theta_2} - \frac{1}{4} e^{\frac{1}{2}\imath \theta_2} - \frac{1}{8} e^{\imath \theta_1} e^{\frac{1}{2} \imath \theta_2}.
\end{align*}}}
The calculations for the restriction operators $R_h^H\big\lvert_v$
and $R_h^H\big\lvert_p$ are similar. In addition, an equivalent calculation (see the work of MacLachlan and Oosterlee \citep{MacLachlanOosterlee2011}) gives the symbols of interpolation operators $\hat{P}_H^h(\bm{\theta})$.
We analyze the two-grid method based on the smoothing operators as introduced in the previous section. Based on the updated restriction and prolongation operators, we perform the analysis of the two-grid method with triad-wise and Vanka smoothing procedures. The asymptotic two-grid convergence factor is given by the spectral radius
\begin{align*}
\rho(E_{TG}):= \sup_{\bm{\theta} \in T^{\text{low}}} \rho\left(\hat{E}_{TG}(\bm{\theta})\right),
\end{align*}
where 
\begin{align*}
\hat{E}_{TG}(\bm{\theta})&=\left(\hat{\mathcal{S}}_h(\bm{\theta})\right)^{\nu_2} \hat{B}_h^H(\bm{\theta})\left(\hat{\mathcal{S}}_h(\bm{\theta})\right)^{\nu_1}\\  &=\left(\hat{\mathcal{S}}_h(\bm{\theta})\right)^{\nu_2} \left(\hat{I}-\hat{P}_H^h(\bm{\theta})\left(\hat{L}_H(2\bm{\theta})\right)^{-1}\hat{R}_h^H(\bm{\theta})\hat{L}_h(\bm{\theta})\right) \left(\hat{\mathcal{S}}_h(\bm{\theta})\right)^{\nu_1}.
\end{align*}

\subsection{Analysis results}
\label{sec:anares}

\begin{figure}[h]
\centering
\subcaptionbox[]{\label{subfig:4a} Vanka smoother} 
[0.3\textwidth] 
{\includegraphics[width=5.5 cm]{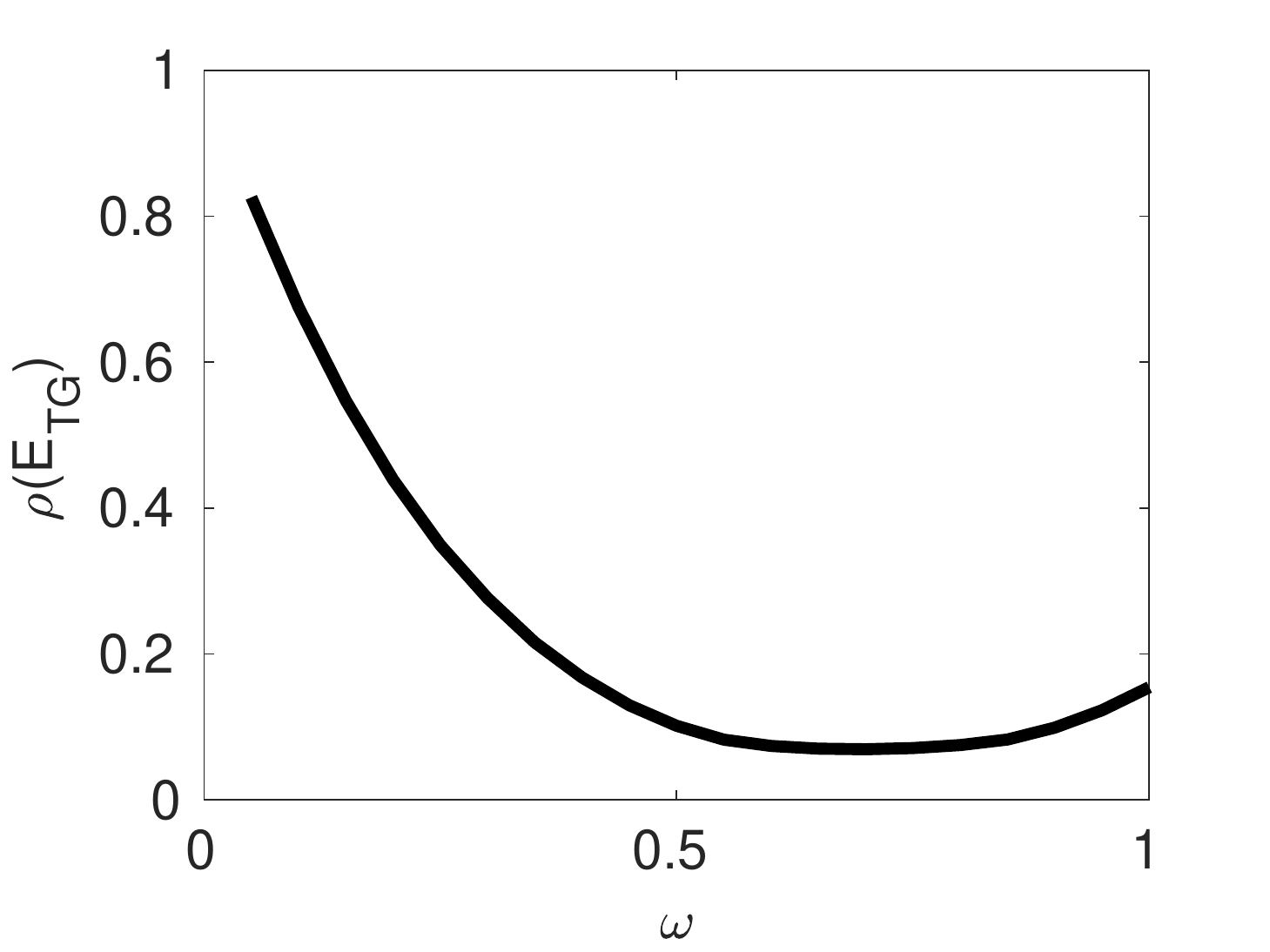}}
\subcaptionbox[]{\label{subfig:4b} triad-Gauss-Seidel smoother} 
[0.3\textwidth] 
{\includegraphics[width=5.5 cm]{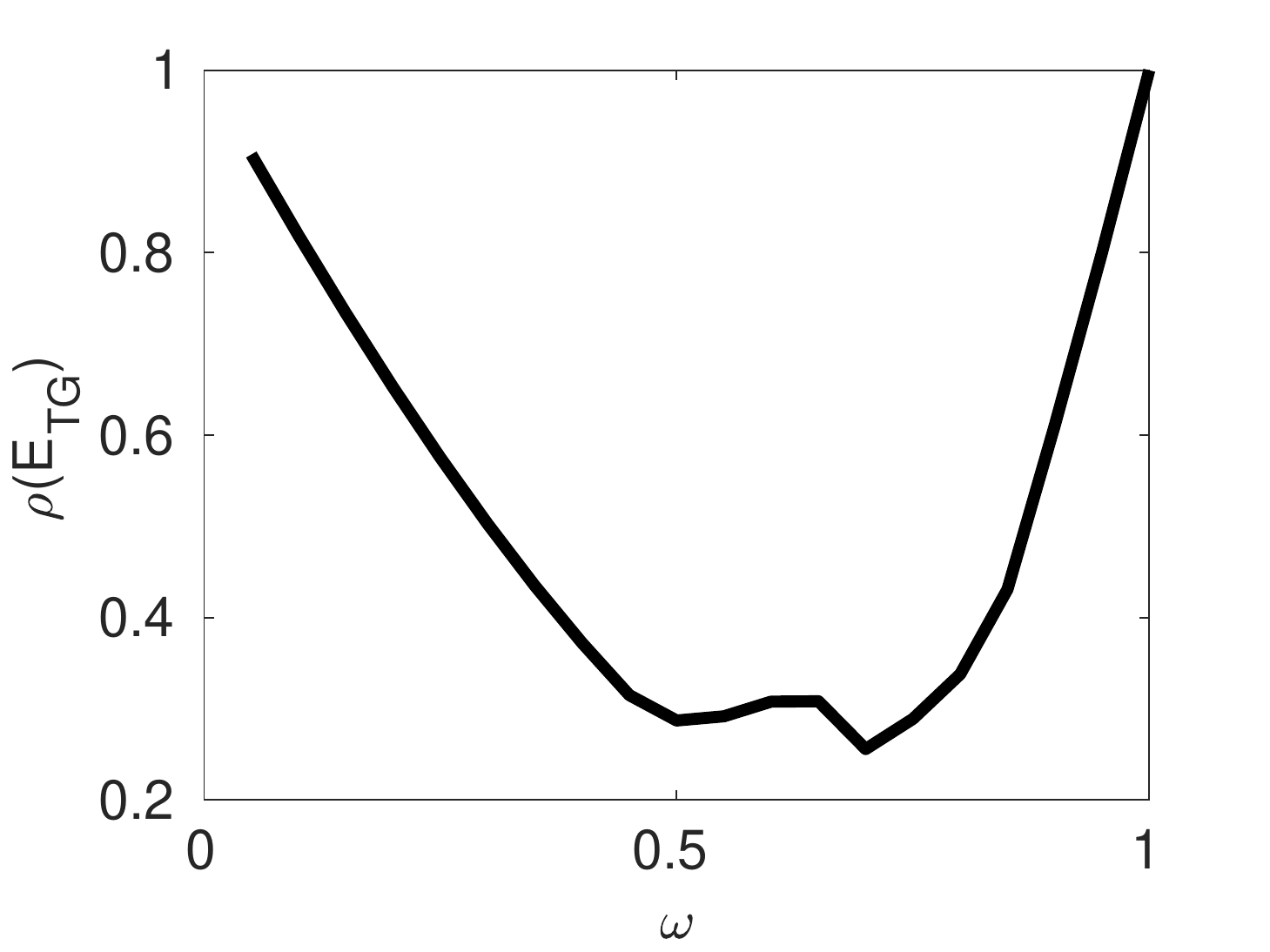}}
\subcaptionbox[]{\label{subfig:4c} triad-Jacobi smoother} 
[0.3\textwidth] 
{\includegraphics[width=5.5 cm]{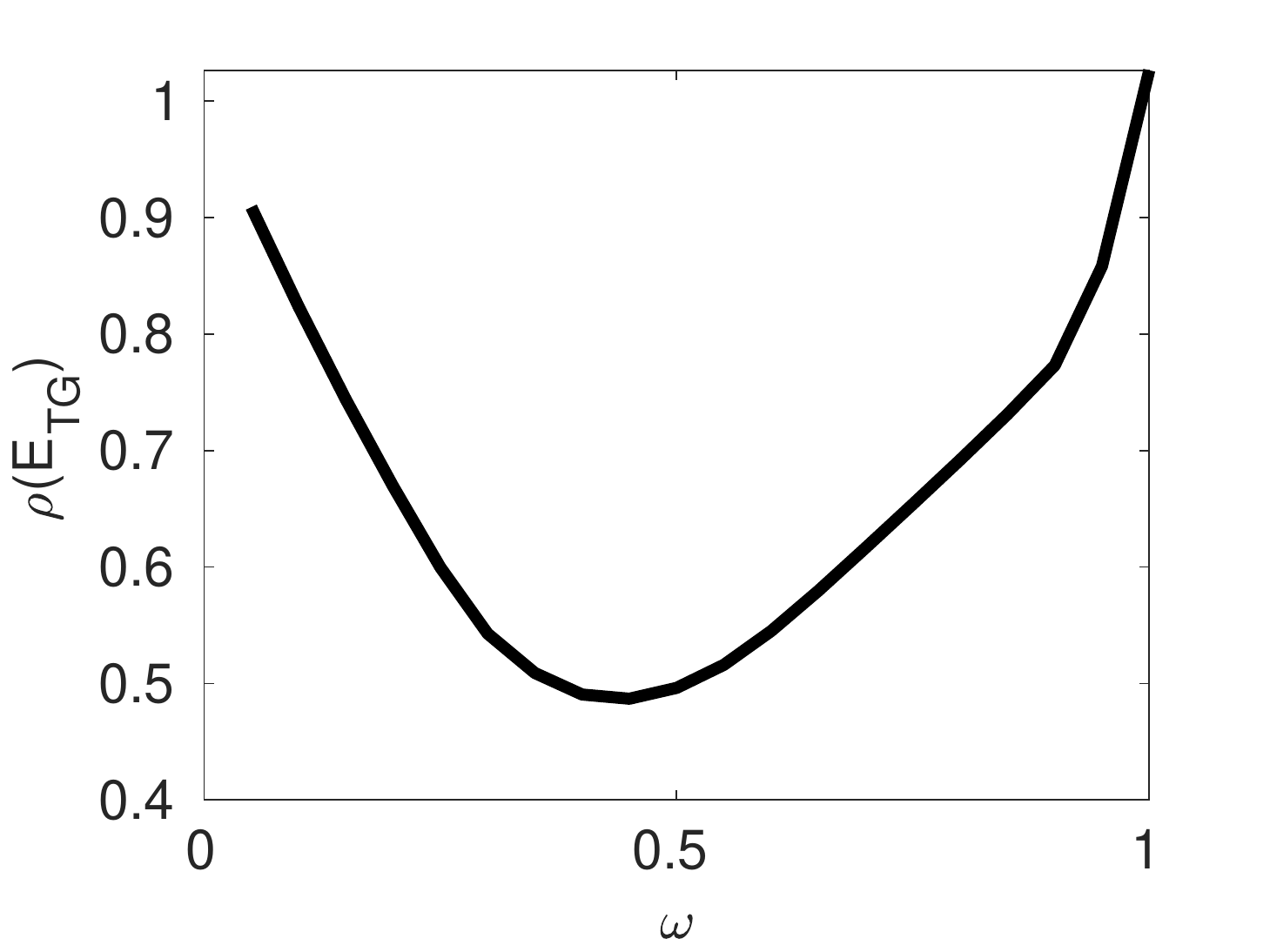}}
\caption{Parameter study of the two-grid method for different smoothers.}
\label{fig:4}
\end{figure}

\begin{figure}[h]
\centering
\subcaptionbox[]{\label{subfig:5a}} 
[0.45\textwidth] 
{\includegraphics[width=5.35 cm]{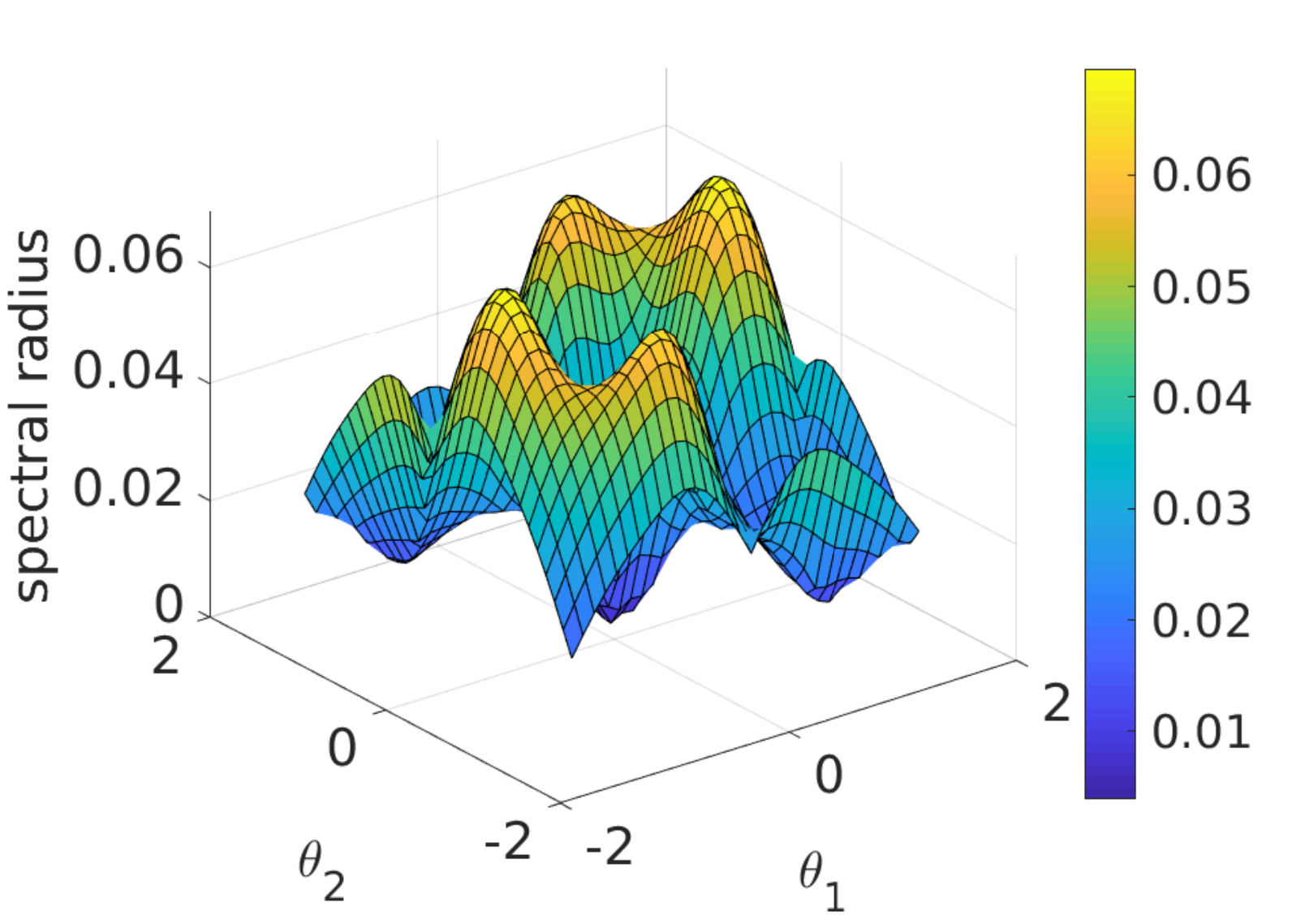}}
\subcaptionbox[]{\label{subfig:5b}} 
[0.45\textwidth] 
{\includegraphics[width=5.35 cm]{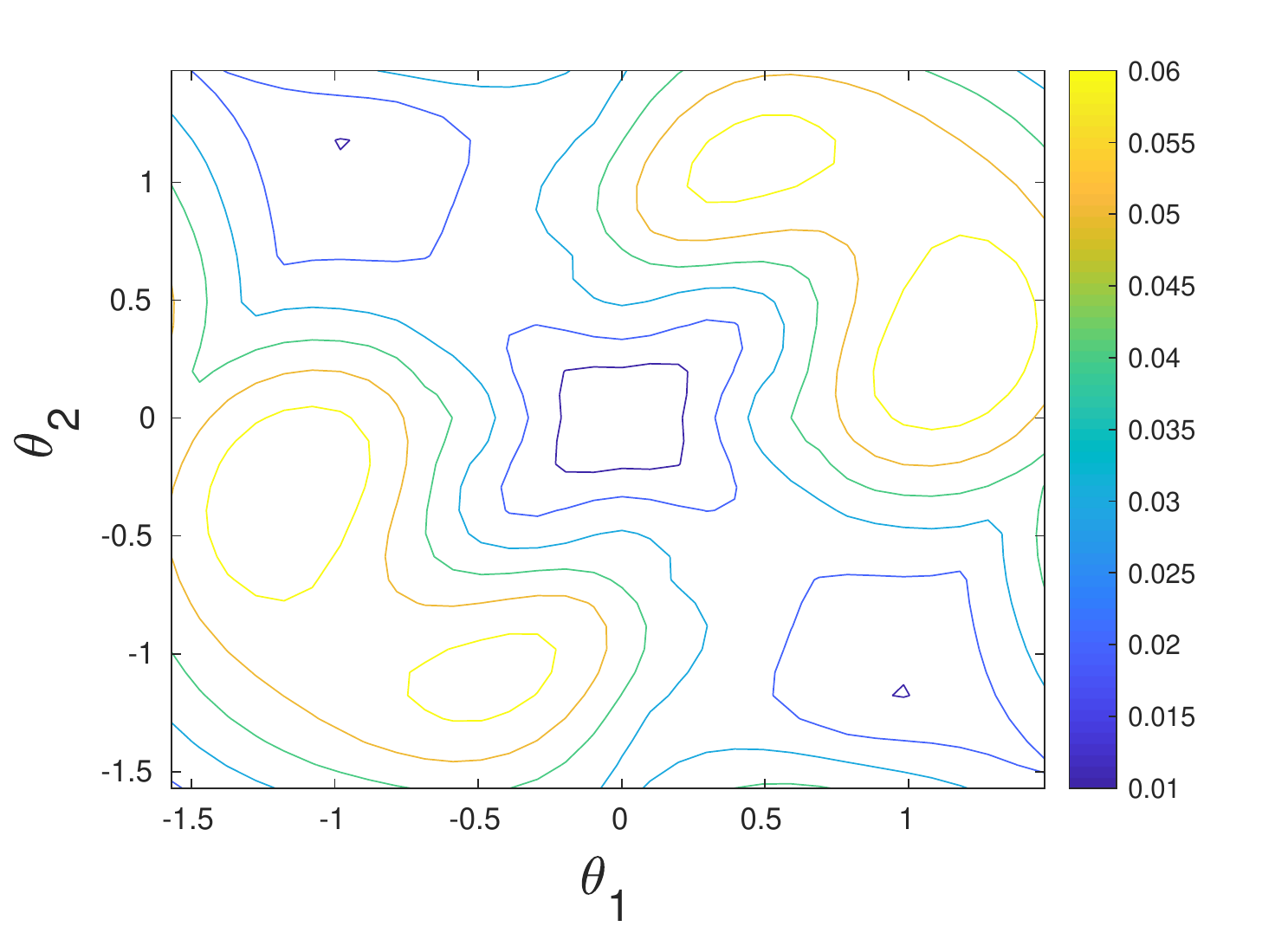}}
\caption{The spectral radius of the two-grid symbol for the Vanka smoother with $\omega = 0.7$ and different frequencies $\bm{\theta}$.}
\label{fig:5}
\end{figure}

\begin{figure}[h]
\centering
\subcaptionbox[]{\label{subfig:6a}} 
[0.45\textwidth] 
{\includegraphics[width=5.35 cm]{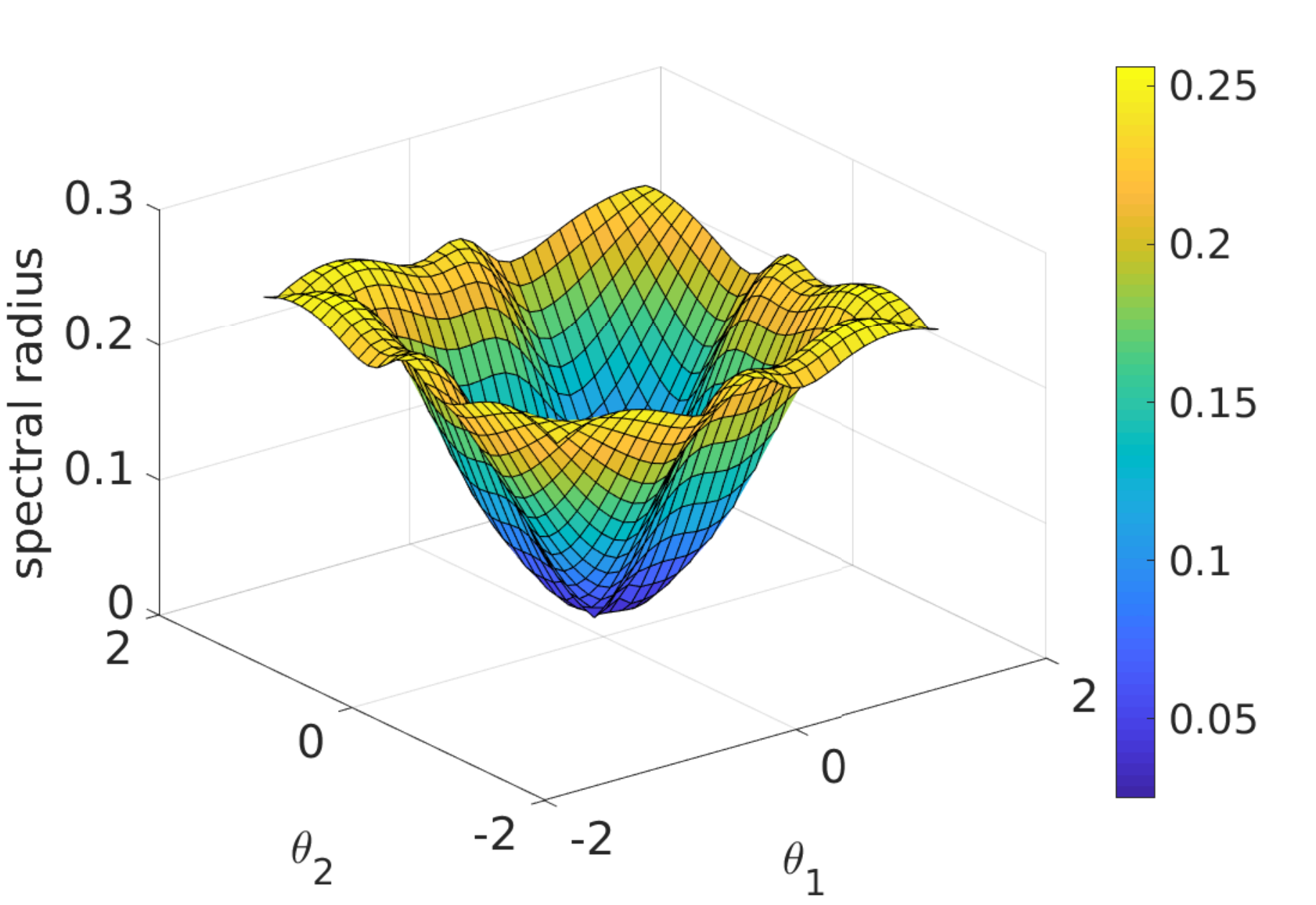}}
\subcaptionbox[]{\label{subfig:6b}} 
[0.45\textwidth] 
{\includegraphics[width=5.35 cm]{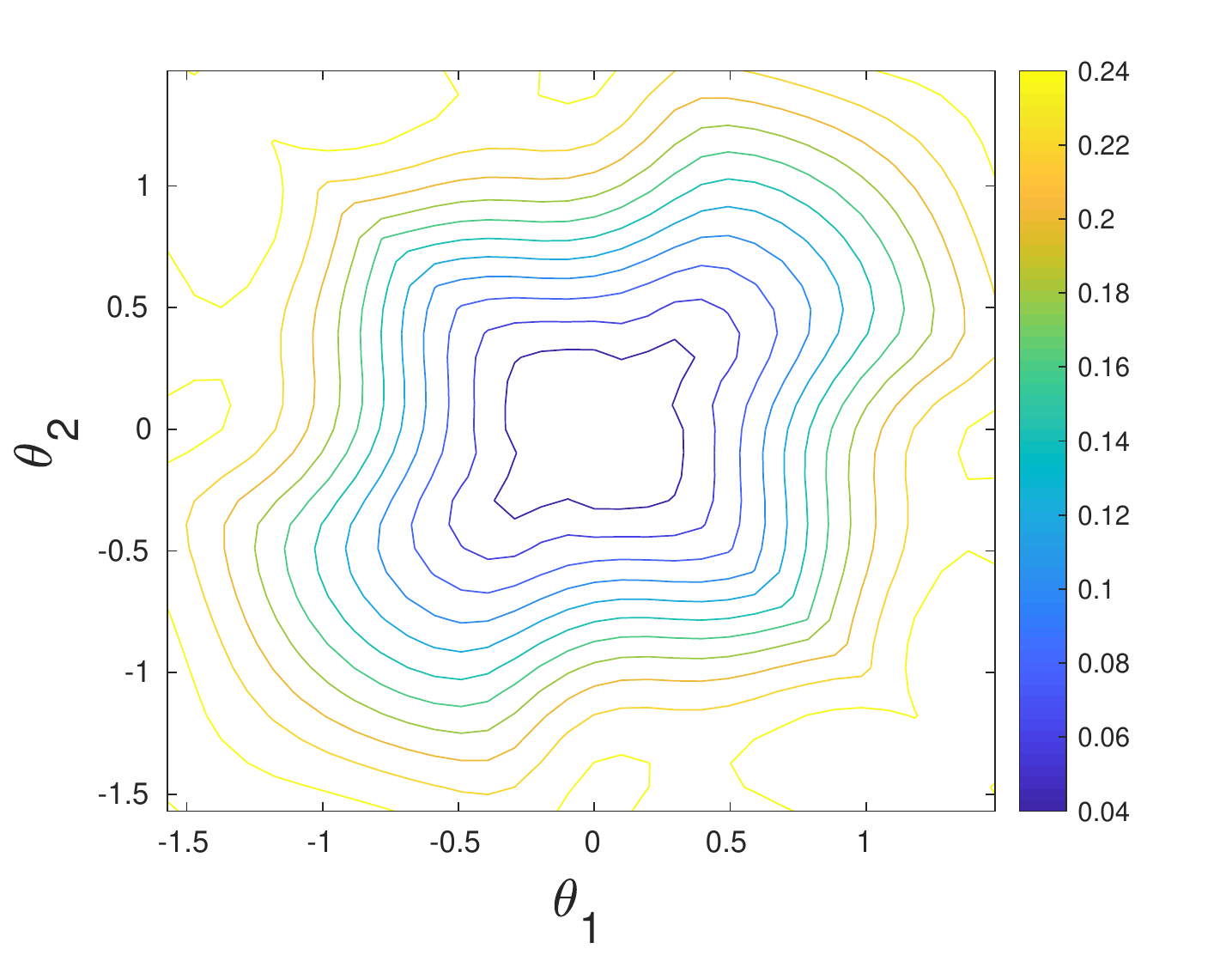}}
\caption{The spectral radius of the two-grid symbol for the triad-wise relaxation based on a weighted Gauss-Seidel updating process with $\omega = 0.7$ and different frequencies $\bm{\theta}$.}
\label{fig:6}
\end{figure}

\begin{figure}[h]
\centering
\subcaptionbox[]{\label{subfig:7a}} 
[0.45\textwidth] 
{\includegraphics[width=5.35 cm]{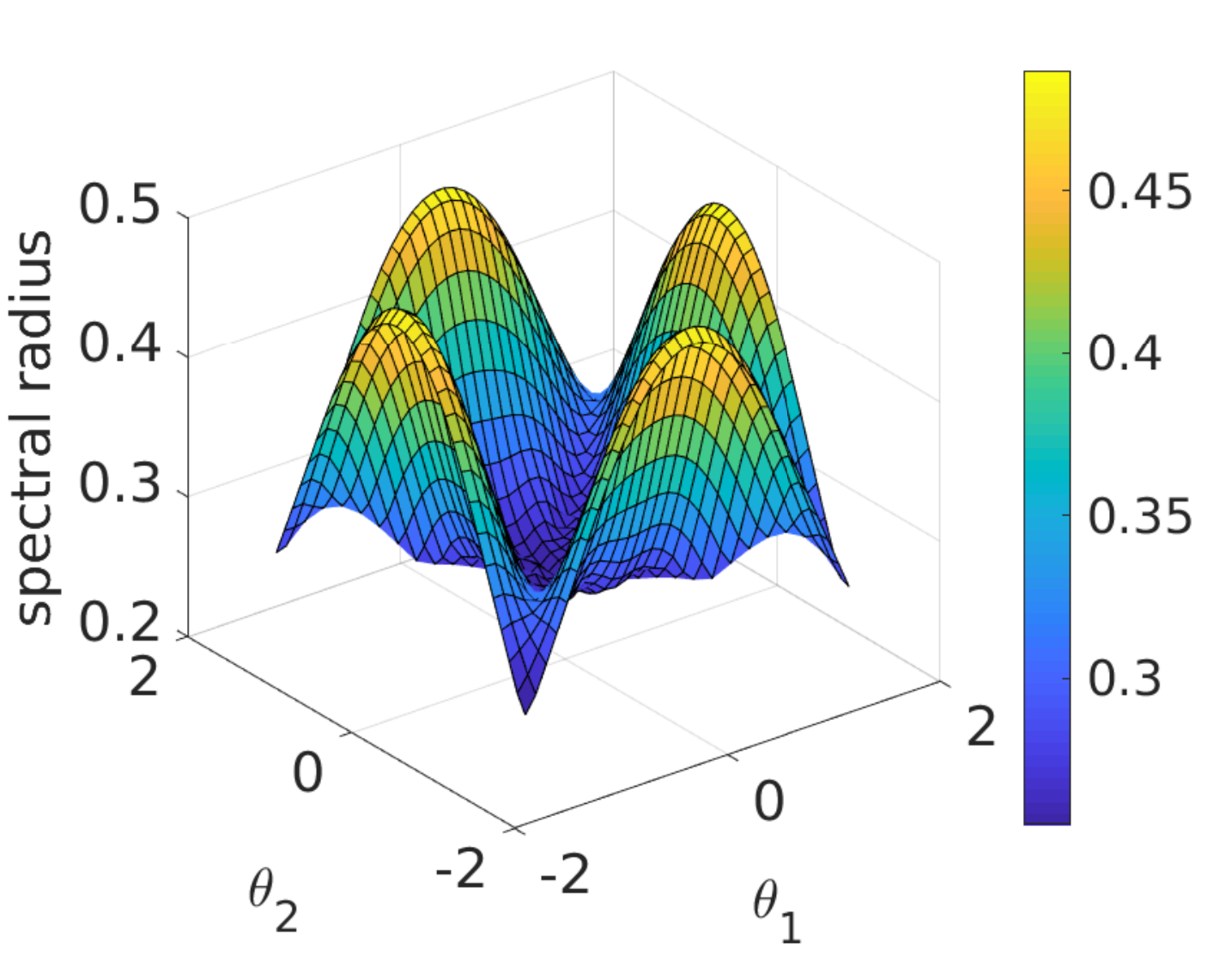}}
\subcaptionbox[]{\label{subfig:7b}} 
[0.45\textwidth] 
{\includegraphics[width=5.35 cm]{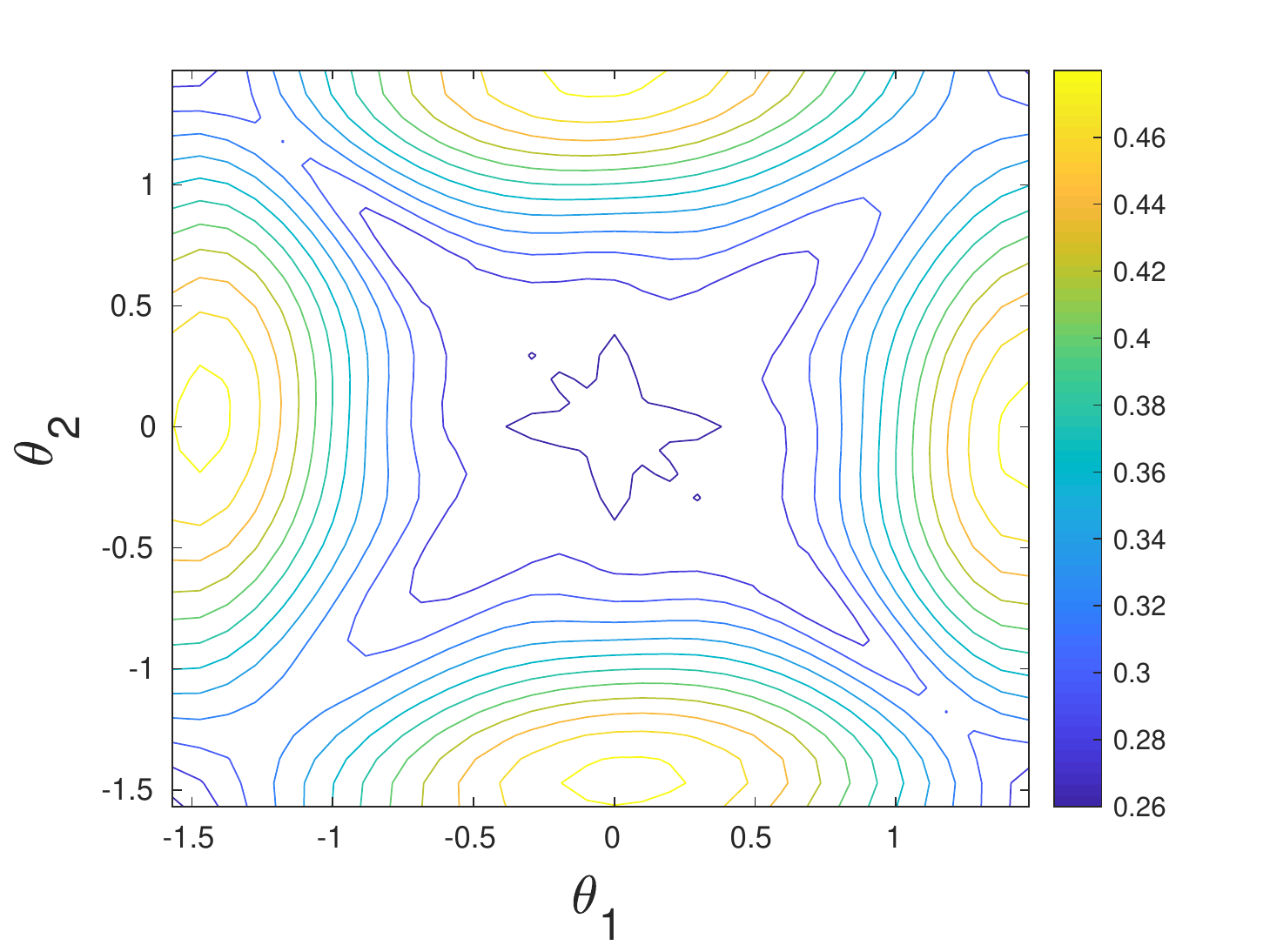}}
\caption{The spectral radius of the two-grid symbol for the triad-wise relaxation based on a weighted Jacobi updating process with $\omega = 0.45$ and different frequencies $\bm{\theta}$.}
\label{fig:7}
\end{figure}

In this section, we present LFA two-grid results with a focus on bilinear interpolation and a Galerkin coarse grid operator. A comparison in \citep{clausphd} shows that the combination of these operators promises best results for the Stokes equations. We compare the results obtained by Vanka and triad-wise smoothers. We use a $33 \times 33$ staggered grid and two pre- and postsmoothing steps ($\nu$=2) during the course of this section. Figure \ref{fig:4} presents the LFA convergence factors for the Vanka and triad-wise smoother as a function of the weighting parameter $\omega$. This shows the sensitivity of performance to parameter choice. \\
On the basis of these parameter studies, c.f.\ Figure \ref{fig:4}, we figured out that the weighting factor $\omega=0.7$ gives the most promising smoothing behavior for the Vanka and triad-Gauss-Seidel(GS) smoother. While for the triad-Jacobi updating process, we stick to the weighting factor $\omega=0.45$ in the following. As visualized in Figures \ref{fig:5}, \ref{fig:6} and \ref{fig:7}, the above mentioned combination of MG components lead to the convergence factor
\begin{align*}
\rho(E_{TG})=0.08,
\end{align*}
for the Vanka smoother,
\begin{align*}
\rho(E_{TG})=0.26,
\end{align*}
for the triad-GS smoother, and
\begin{align*}
\rho(E_{TG})=0.49,
\end{align*}
for the triad-Jacobi smoother. 

The performance predicted by LFA shows the potential of the triad-wise smoother.

\section{Numerical examples}\label{sec4}

\begin{figure}[h]
\centering
\includegraphics[width=13 cm]{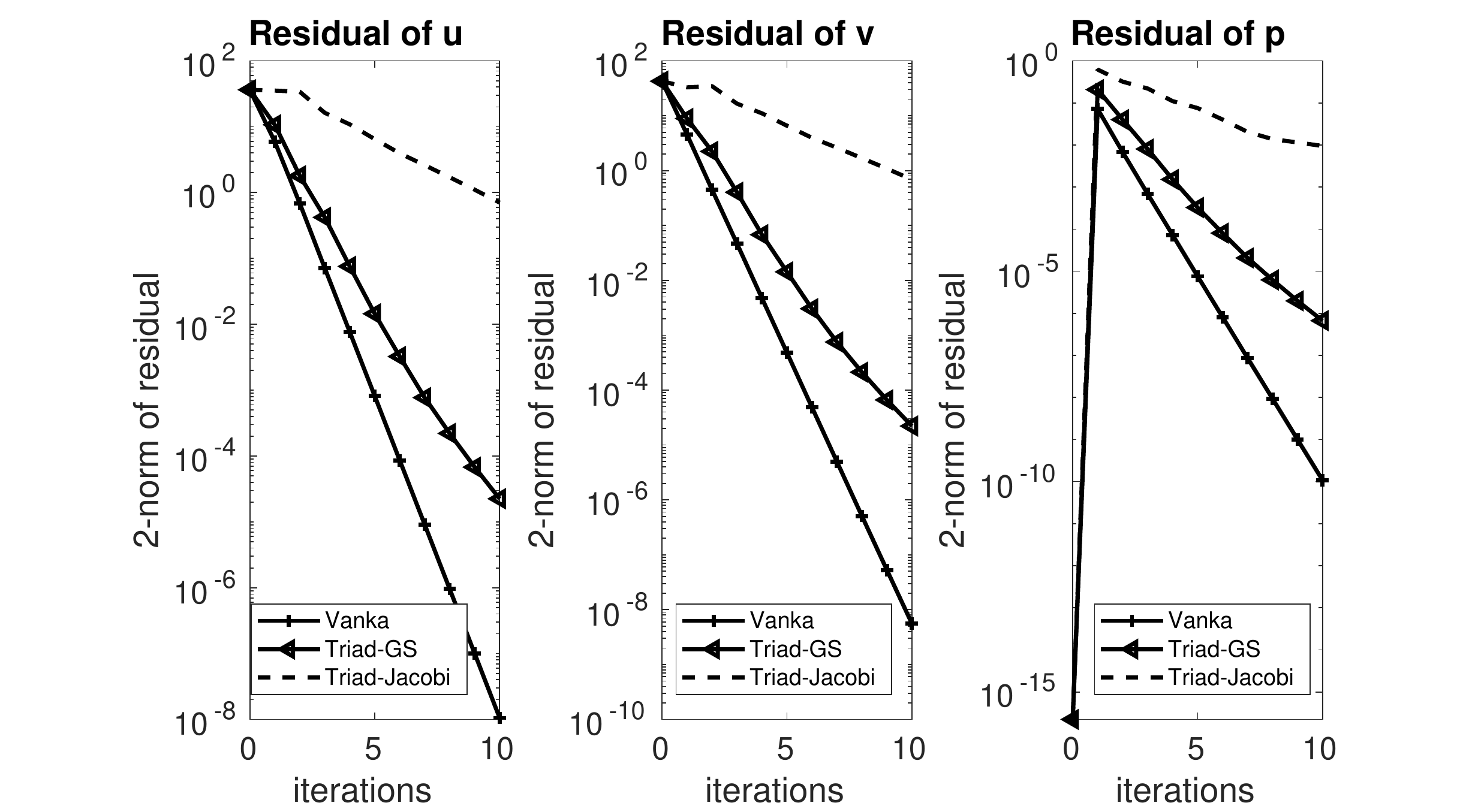}
\caption{Convergence behavior of the multigrid method for the Stokes system with periodic boundary conditions.}
\label{fig:8}
\end{figure}

In this section, we first present convergence results of a two-grid and a V-cycle multigrid method for the Stokes equations. We start with periodic boundary conditions and observe convergence behavior as predicted by LFA in the previous section. 
A two-grid method for a homogeneous problem with periodic boundary conditions based on a $33 \times 33$ staggered grid is employed. We apply 20 two-grid cycles with two pre- and postsmoothing steps ($\nu :=\nu_1=\nu_2=2$) and start with a random initial guess. We measure the convergence via the 2-norm of the error, i.e.\ we compute $\frac{||\mathbf{e}^k||_2}{||\mathbf{e}^{k-1}||_2}$ after $k=20$ two-grid cycles. This gives a convergence factor of $\rho=0.08$ for the Vanka smoother with $\omega = 0.7$, $\rho=0.24$ for the weighted ($\omega = 0.7$) triad-GS smoother and $\rho=0.43$ for the weighted ($\omega = 0.45$) triad-Jacobi smoother. The numerical results confirm the LFA convergence results as given in Section \ref{sec:anares}. In addition, Figure \ref{fig:8} shows convergence results for the Vanka and the triad-wise smoothers for the multigrid method.

Figure \ref{fig:8} is based on a $33 \times 33$ staggered fine grid and a $3 \times 3$ coarse grid on the domain $\Omega=(0,1)^2$ discretized using the MAC scheme. We use two pre- and postsmoothing steps ($\nu=2$) and a zero initial guess. The right-hand side is set to 
\begin{align*}
\mathbf{f}_1 & = 8\pi^2\cdot\sin(2 \pi {\mathrm{x}^1_1})\cdot\sin(2 \pi {\mathrm{x}^1_2})
		 -2\pi \cdot\sin(2 \pi {\mathrm{x}^1_1})\cdot\sin(2 \pi {\mathrm{x}^1_2}),\\
    \mathbf{f}_2 & = 8\pi^2\cdot\cos(2 \pi {\mathrm{x}_1^2})\cdot\cos(2 \pi {\mathrm{x}^2_2})
		+2\pi\cdot\cos(2 \pi {\mathrm{x}^2_1})\cdot\cos(2 \pi {\mathrm{x}^2_2}),\\
\end{align*}
where the values $\mathrm{x}^j$ correspond to the grid points of the discretized domain\linebreak $\Omega_h^j =\{ \bm{x}^j := \bm{k} \mathbf{h} + \delta^j, \bm{k} \in \mathbb{Z}^2 \}$, with $\delta^1 = (0, h/2)$ and $\delta^2 = (h/2, 0)$.
 Figure \ref{fig:8} presents the convergence factors corresponding to the different unknowns $u,v,p$ individually. The convergence factors are obtained by computing the 2-norm of the residual for different numbers of multigrid V-cycles (iterations) for the Stokes equations, i.e.\ we compute $||\mathbf{r}^k||_2$ after each V-cycle. We notice good convergence factors for the Vanka and the triad-GS relaxation methods. However, the Vanka smoother converges even faster than the triad-GS smoother. The triad-Jacobi procedure does not lead to a good convergence rate.\\ 
Next, we show convergence results for the Stokes system with Dirichlet boundary conditions. For the triad-wise smoother we focus on a Gauss-Seidel updating process, since we have seen that the triad-Jacobi smoother does not lead to good convergence results. For Dirichlet boundary conditions triad-Jacobi downgrades to a divergent method. Again, practical two-grid convergence factors for a homogeneous problem discretized on a $33 \times 33$ grid are given. We apply 20 two-grid cycles, set $\nu=2$ and start with a random initial guess. That gives a convergence of $\rho=0.36$ for the triad-wise smoother and $\rho=0.10$ for the Vanka smoother. Again, the Vanka smoother leads to good convergence results. In contrast, the triad smoother does not show satisfactory convergence for Dirichlet boundary conditions. For this reason, we investigate different possibilities to adapt the smoother. We start with additional smoothing steps and vary the order of updating different blocks. We use a lexicographical updating process, i.e. forward direction, and a backward direction as well as red-black(RB)-coloring. The results can be found in Table \ref{table:tab1}. 

\begin{table}[h]
\centering
\renewcommand{\arraystretch}{2}
\begin{tabular}{ c|c c c c }
 &  $\omega$-GS (forward) & $\omega$-GS (backward) & $\omega$-GS-RB (backward) & $\omega$-GS, $\nu=6$ (backward) \\
 \hline
$\mu$ & 0.62 & 0.36 & 0.58 & 0.29 \\

\end{tabular}
\caption{Convergence factors for different triad-type relaxation methods with Dirichlet boundary conditions.}\label{table:tab1}
\end{table}

\begin{figure}[h]
\centering
\includegraphics[width=3.3 cm]{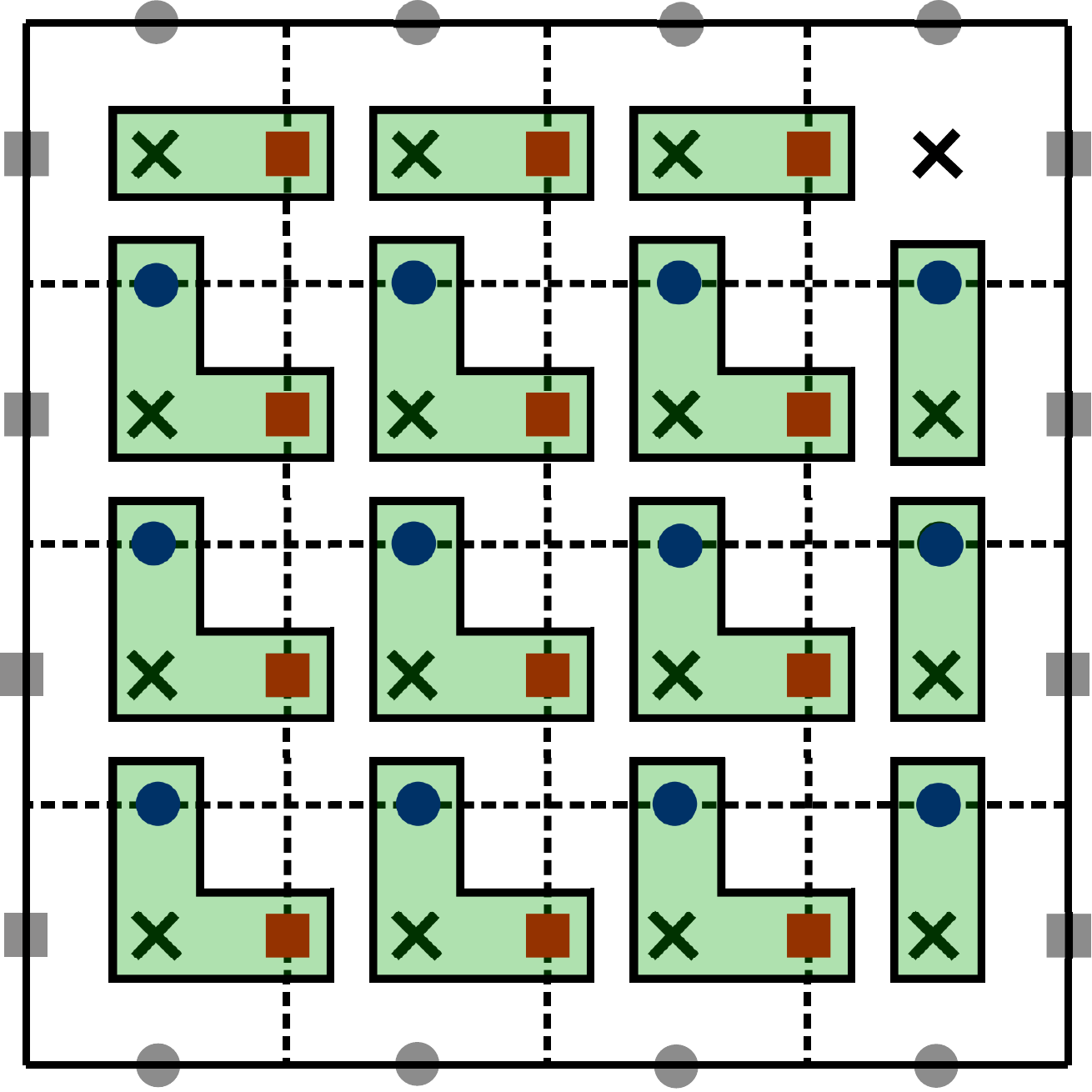}
\caption{Examplary choice of triad blocks with Dirichlet boundary conditions.}
\label{fig:9}
\end{figure}

Table \ref{table:tab1} shows that a different ordering results in small effects only. An increase in the number of smoothing steps leads to better convergence results. However, we obtain a convergence factor of $\rho=0.29$ if we apply six pre- and postsmoothing steps. This is as expensive as using the Vanka smoother but less efficient in terms of the convergence rate.
The poor convergence behavior of the triad-wise smoother results from the boundary treatment, as it does not treat the unknowns at the boundary (or near the boundary) in an appropriate way. Since the unknowns at the boundary are predefined, they are not treated by the relaxation. The idea of the triad-wise smoother is to relax three unknowns at the same time, which is not possible with Dirichlet boundary conditions regardless of the choice of the three unknowns. An illustration with an exemplary choice of triad blocks is given in Figure \ref{fig:9}.

\begin{figure}[h]
\centering
\subcaptionbox[]{\label{subfig:10a}} 
[0.3\textwidth] 
{\includegraphics[width=4 cm]{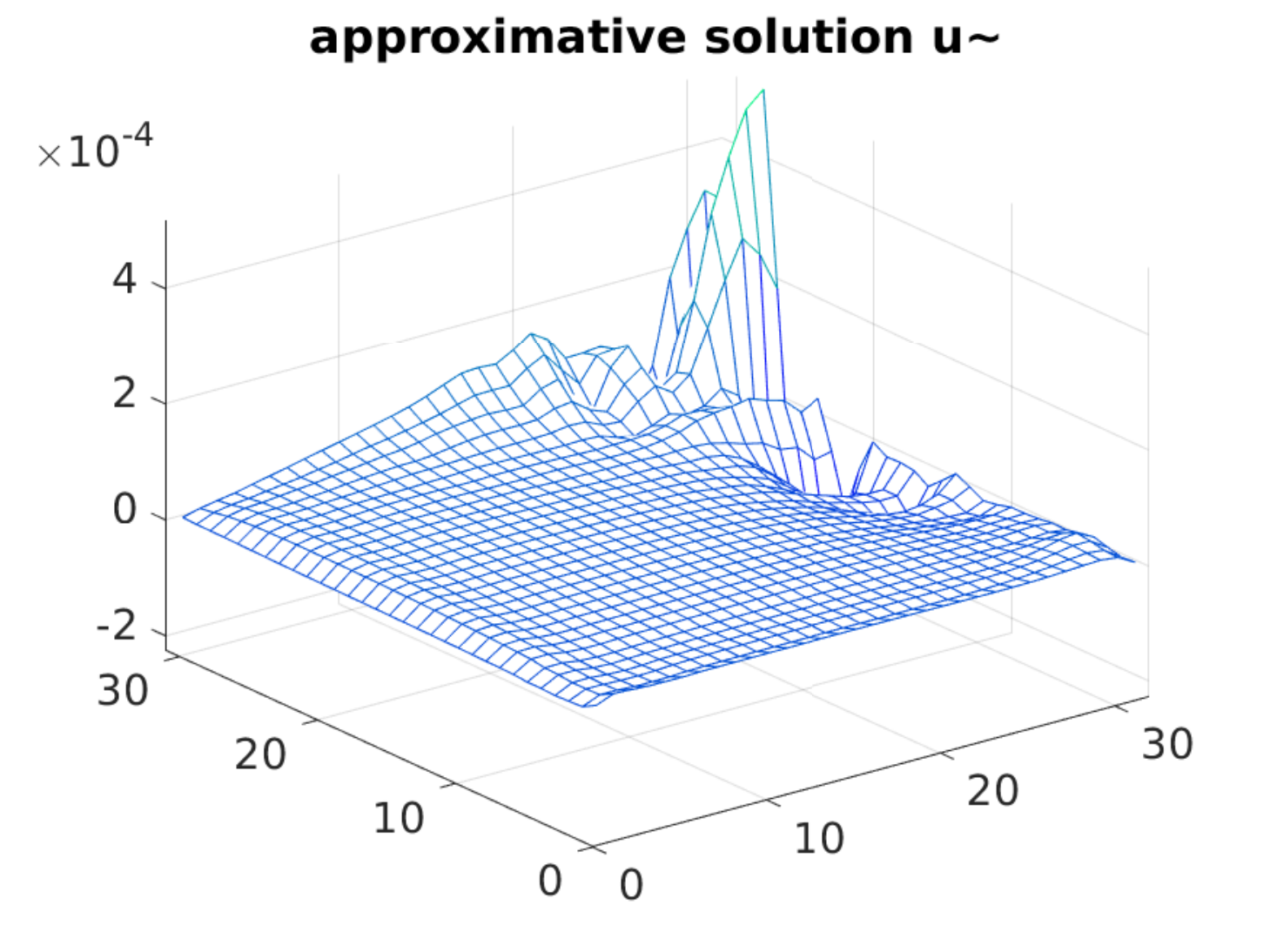}}
\subcaptionbox[]{\label{subfig:10b}} 
[0.3\textwidth] 
{\includegraphics[width=4 cm]{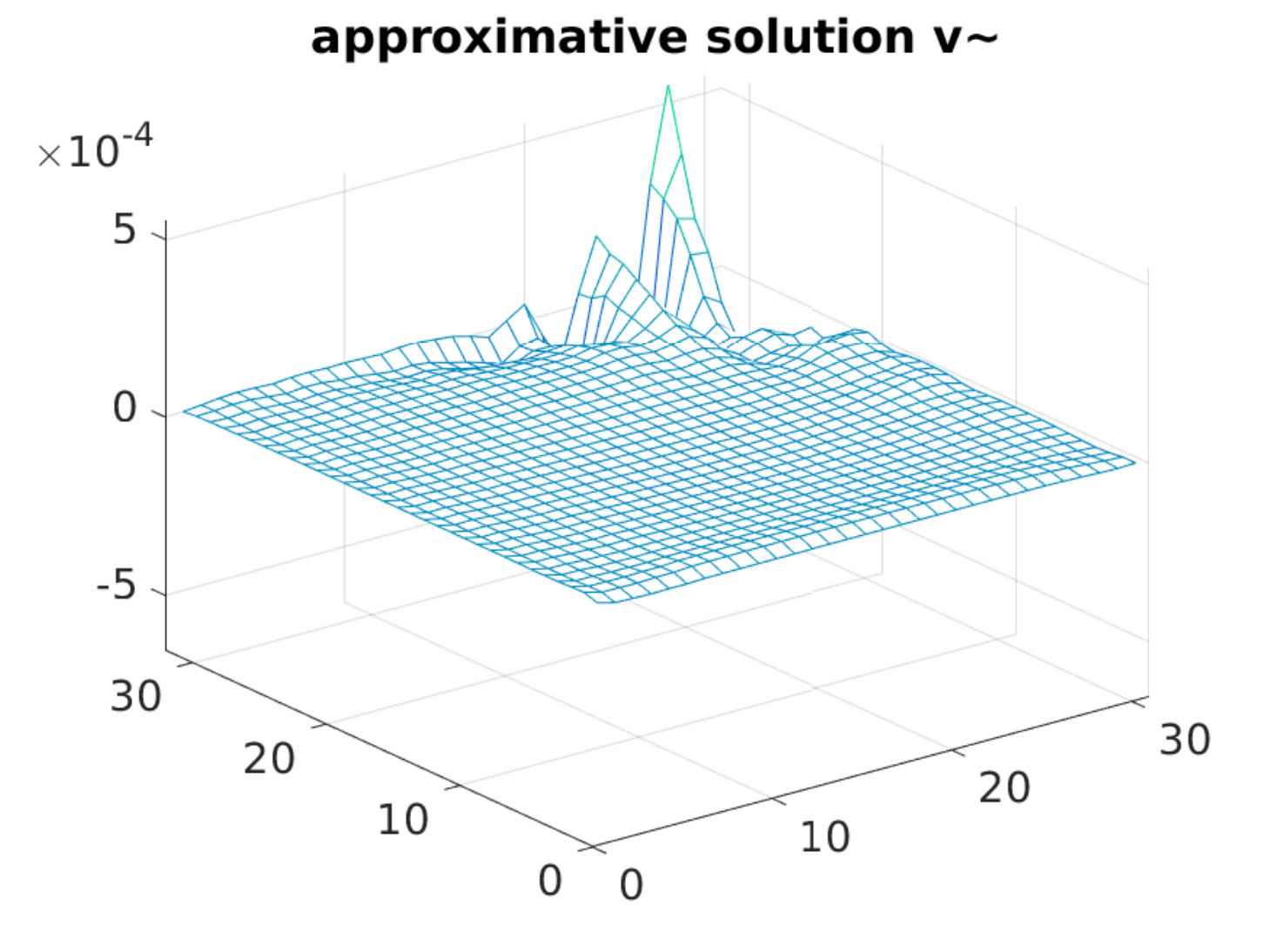}}
\subcaptionbox[]{\label{subfig:10c}} 
[0.3\textwidth] 
{\includegraphics[width=4 cm]{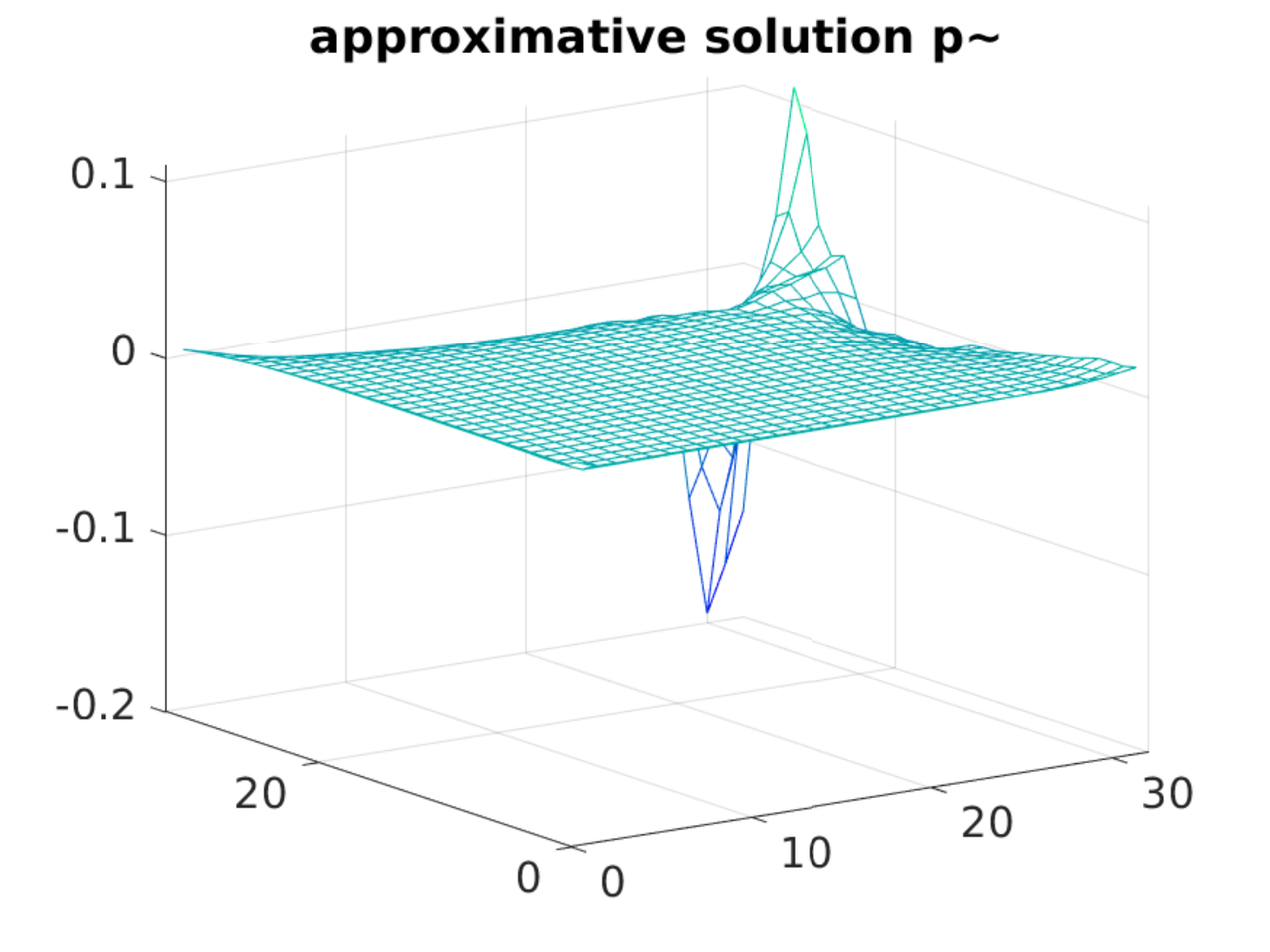}}
\caption{Approximate solution of the Stokes equations with Dirichlet boundary conditions after 20 multigrid cycles based on the triad smoother.}
\label{fig:10}
\end{figure}

We notice the blocks include three unknowns whenever possible. Close to the boundary we have blocks with two unknowns only. At one corner we have one unknown only. Figure \ref{fig:10} illustrates the incorrect treatment by a visualization of the approximate solution after 20 multigrid cycles.

\begin{figure}[h]
\centering
\subcaptionbox[]{\label{subfig:11a}} 
[0.3\textwidth] 
{\includegraphics[width=4 cm]{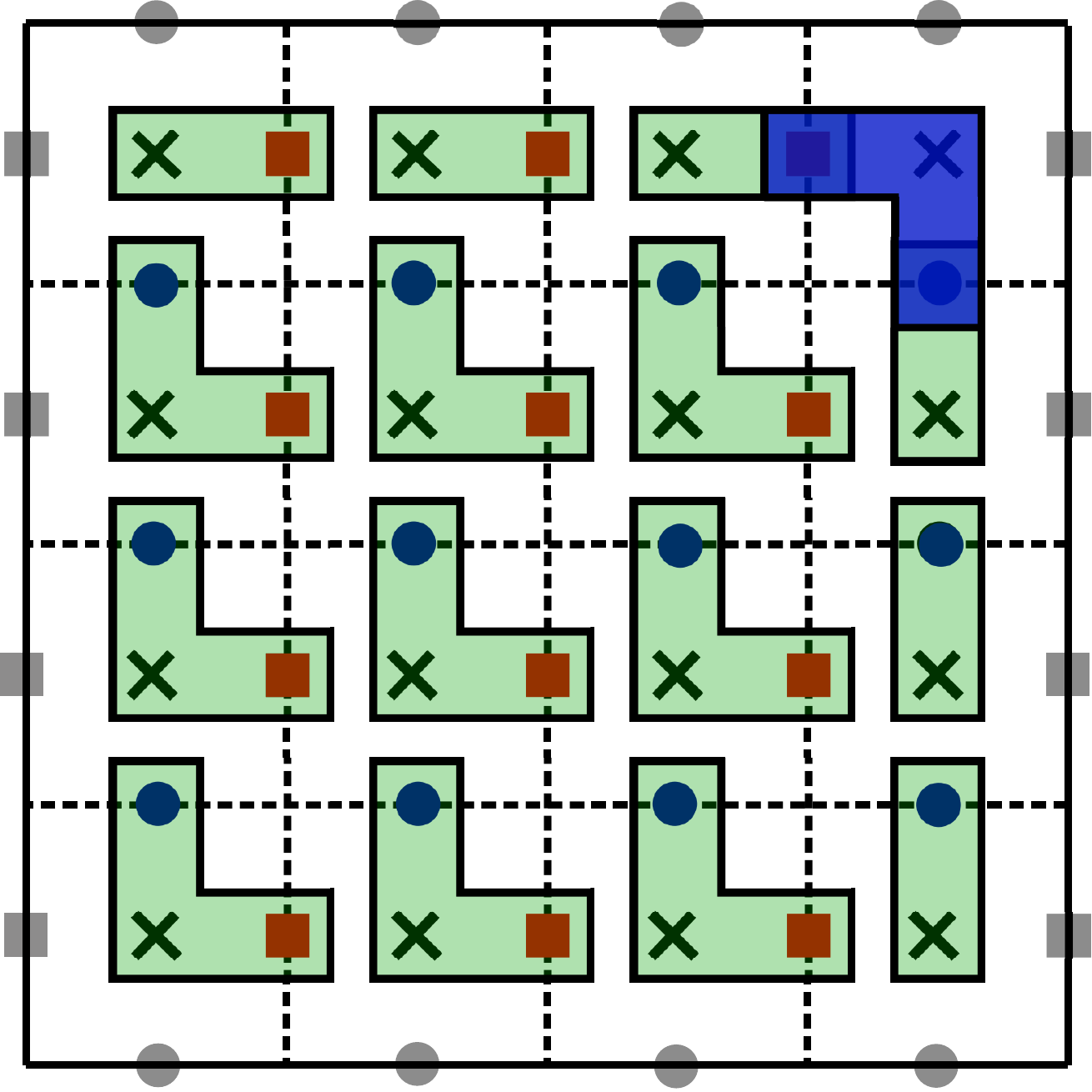}}
\subcaptionbox[]{\label{subfig:11b}} 
[0.3\textwidth] 
{\includegraphics[width=4 cm]{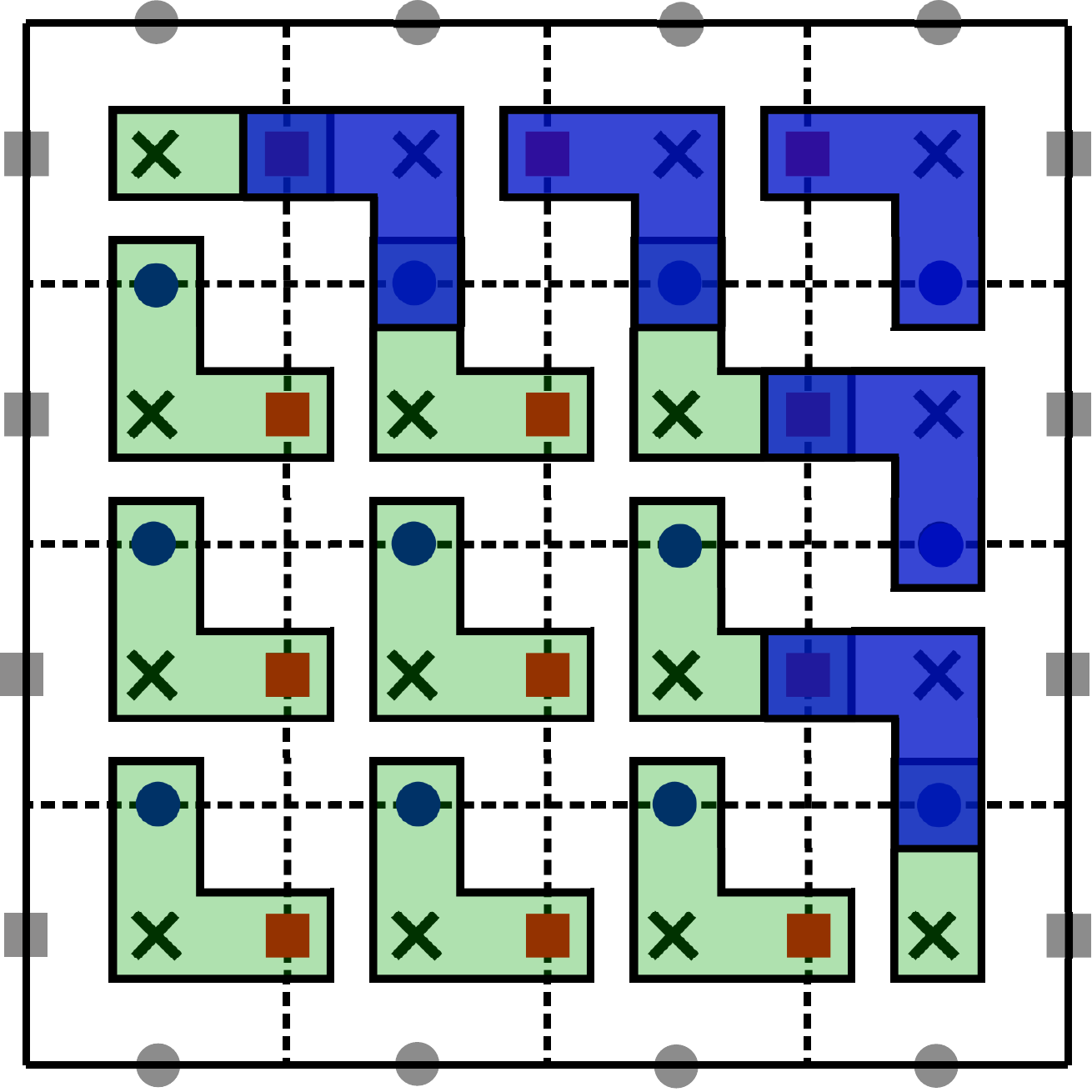}}
\subcaptionbox[]{\label{subfig:11c}} 
[0.3\textwidth] 
{\includegraphics[width=4 cm]{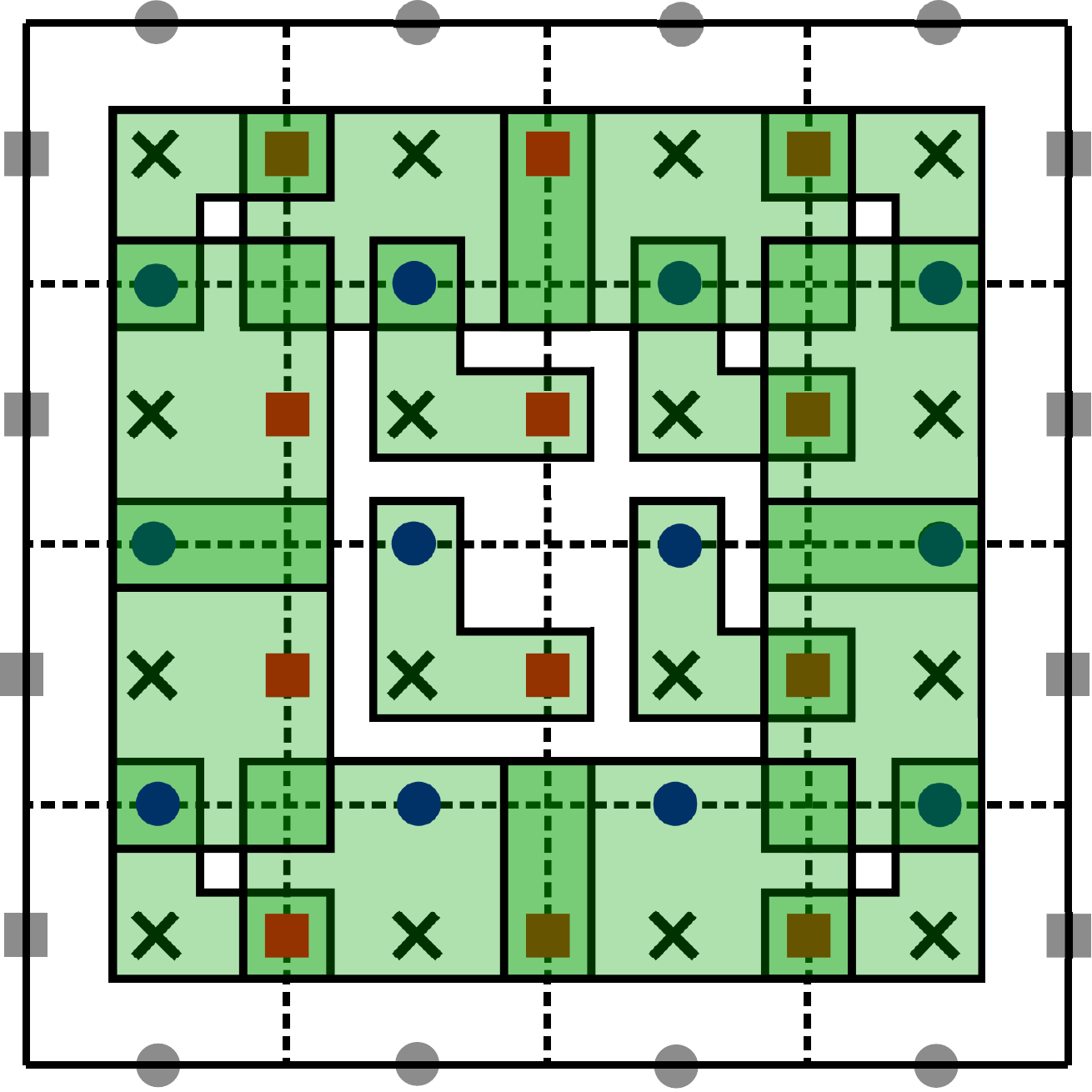}}
\caption{Triad updating blocks with additional overlapping blocks.}
\label{fig:11}
\end{figure}

This observation leads to the following three adaptations of the boundary treatment that are able to improve the smoother but still don't lead to satisfying results. First, after each triad-wise-smoothing cycle, we add an additional smoothing at the critical corner by using an overlapping block as visualized in Figure \ref{subfig:11a}. This gives a convergence factor of $\varphi=0.41$. By adding overlapping blocks not only at the corner but also at the boundaries, see Figure \ref{subfig:11b}, we get a convergence factor of $\varphi=0.37$. By combining the triad updating process at the interior of the grid with the Vanka smoother at the boundary cells, see Figure \ref{subfig:11c}, we improve the convergence factor to $\varphi=0.22$. 

\begin{figure}[h]
\centering
\includegraphics[width=13 cm]{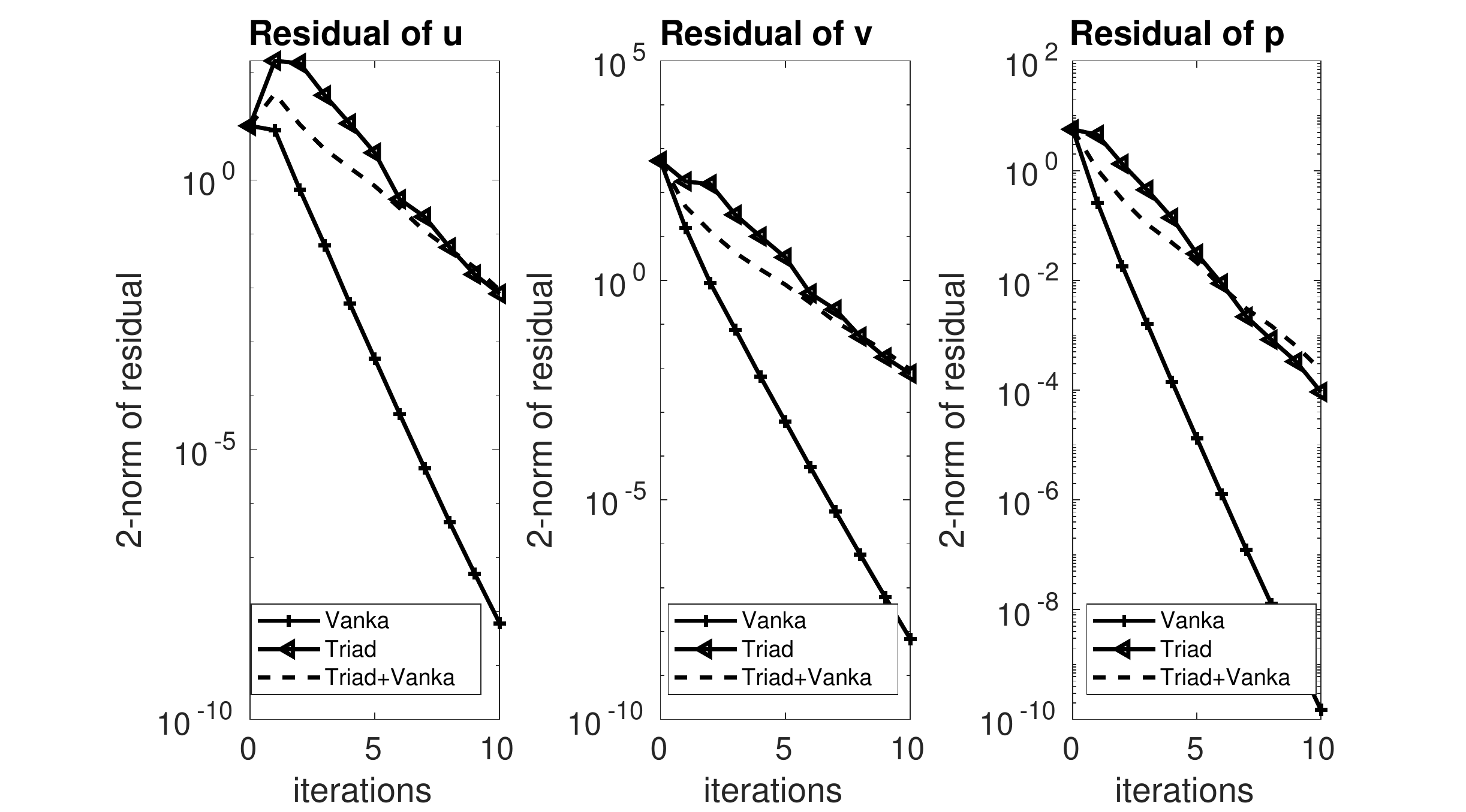}
\caption{Convergence behavior of the multigrid method for the Stokes system with Dirichlet boundary conditions.}
\label{fig:12}
\end{figure}

In addition, we investigate results for a multigrid method as depicted in Figure \ref{fig:12}. We focus on the Vanka smoother and two exemplary triad-wise smoothers. We discretize on a $33 \times 33$ grid with the domain $\Omega=(0,1)^2$, use two pre- and postsmoothing steps ($\nu=2$) and employ the zero solution as an initial guess. 

The right-hand side is set to 
\begin{align*}
    \mathbf{f}_1 & = 2\pi^2\cdot\sin(\pi {\mathrm{x}_1^1})\cdot\sin(\pi {\mathrm{x}^1_2})
         + \pi\cdot\cos(\pi {\mathrm{x}^1_1}),\\
    \mathbf{f}_2 &= 2\pi^2\cdot\cos(\pi {\mathrm{x}^2_1})\cdot\cos(\pi {\mathrm{x}^2_2})
	  - \pi\cdot\sin(\pi {\mathrm{x}^2_2}).\\
\end{align*}
We choose zero boundary conditions for the $x$-component $u$ and $\pm\cos(\pi y)$ and $\pm\cos(\pi x)$ for the $y$-component $v$ of the velocity unknown. Figure \ref{fig:12} shows good convergence for the Vanka smoother (similar to that in the periodic case). In contrast, non of the triad-wise smoothing methods show satisfactory convergence for Dirichlet boundary conditions.

\section{Algorithm development}
\label{sec5}

\begin{figure}[h]
\centering
\subcaptionbox[]{\label{subfig:13a}} 
[0.22\textwidth] 
{\includegraphics[width=0.2\textwidth]{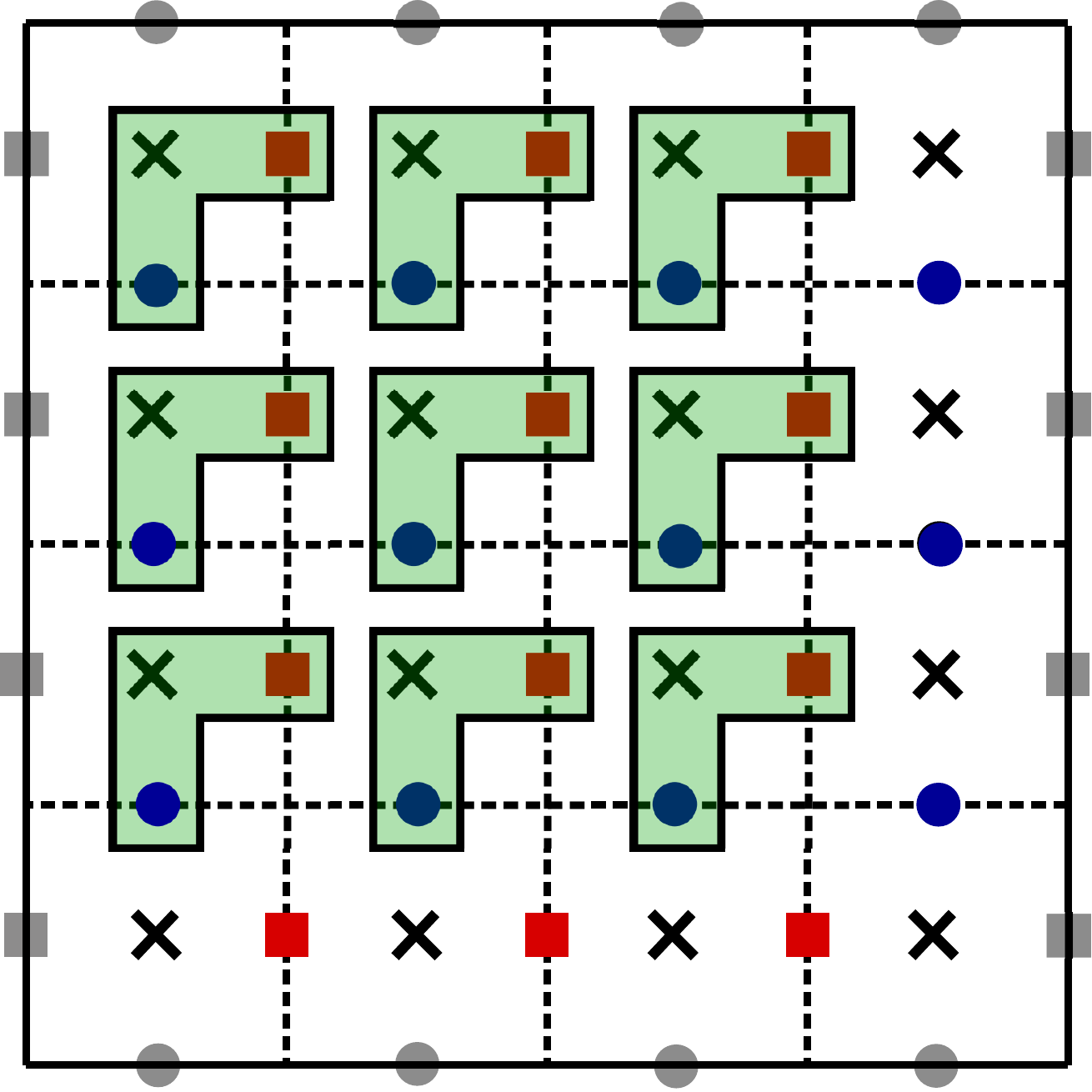}}
\subcaptionbox[]{\label{subfig:13b}} 
[0.22\textwidth] 
{\includegraphics[width=0.2\textwidth]{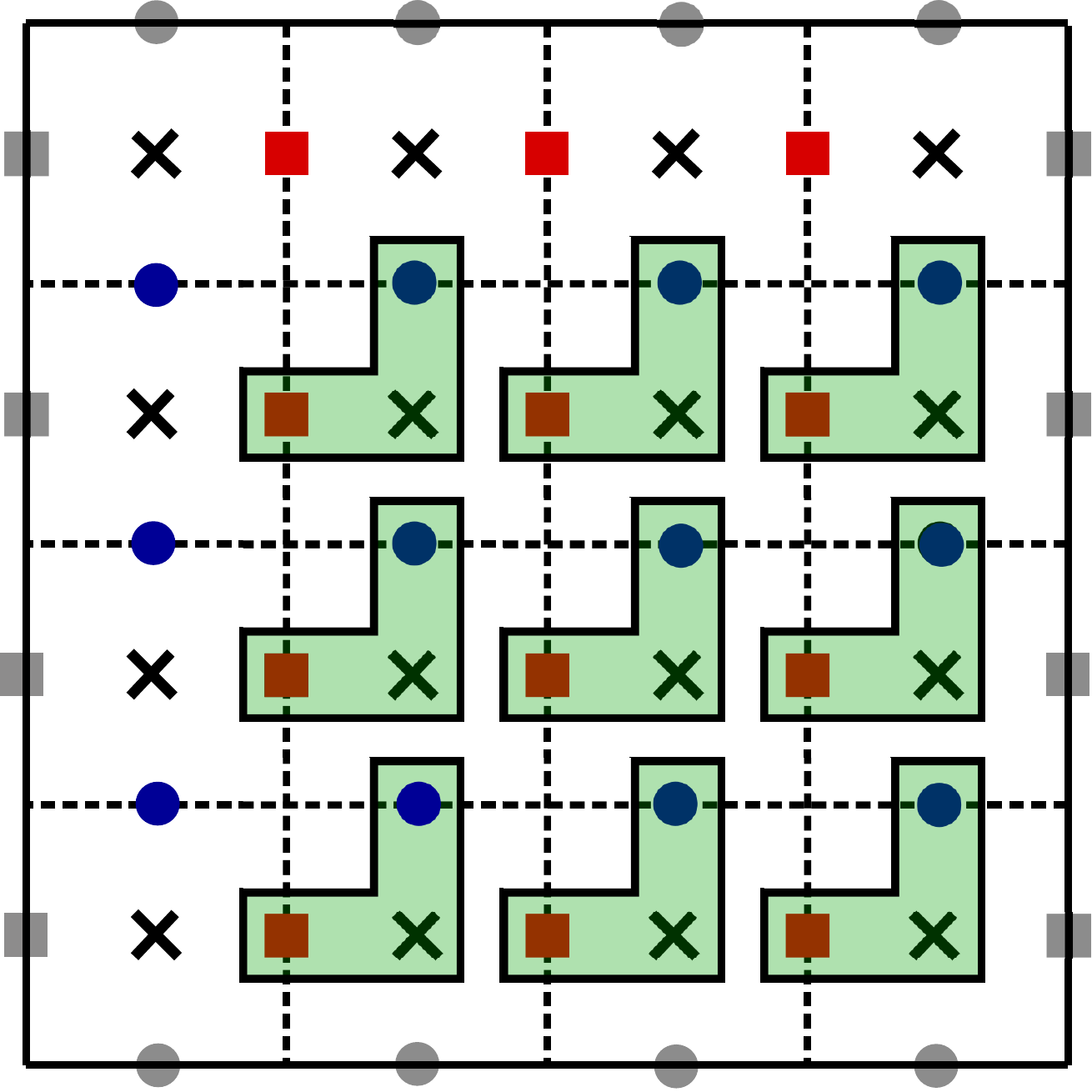} }
\subcaptionbox[]{\label{subfig:13c}} 
[0.22\textwidth] 
{\includegraphics[width=0.2\textwidth]{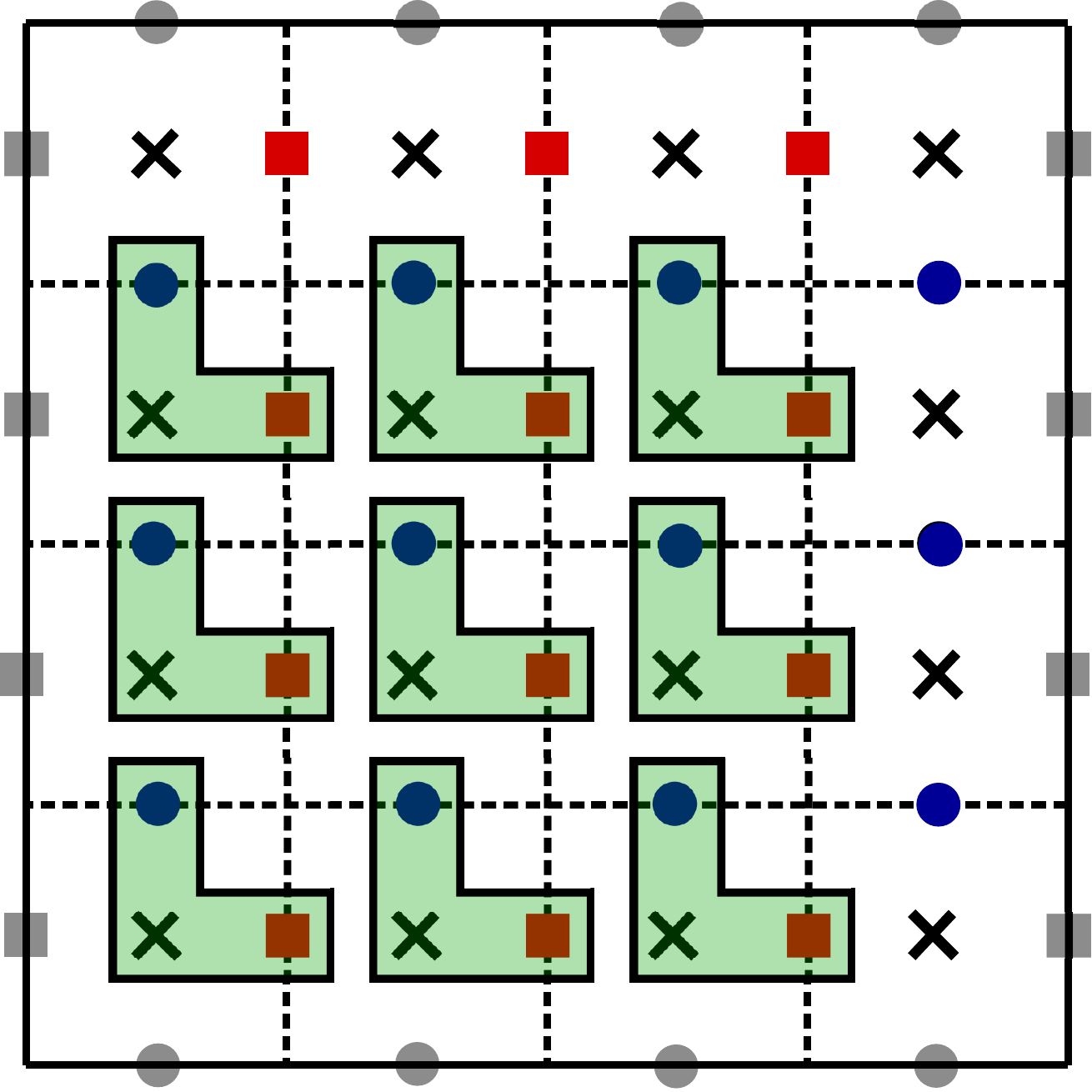}}
\subcaptionbox[]{\label{subfig:13d}} 
[0.22\textwidth] 
{\includegraphics[width=0.2\textwidth]{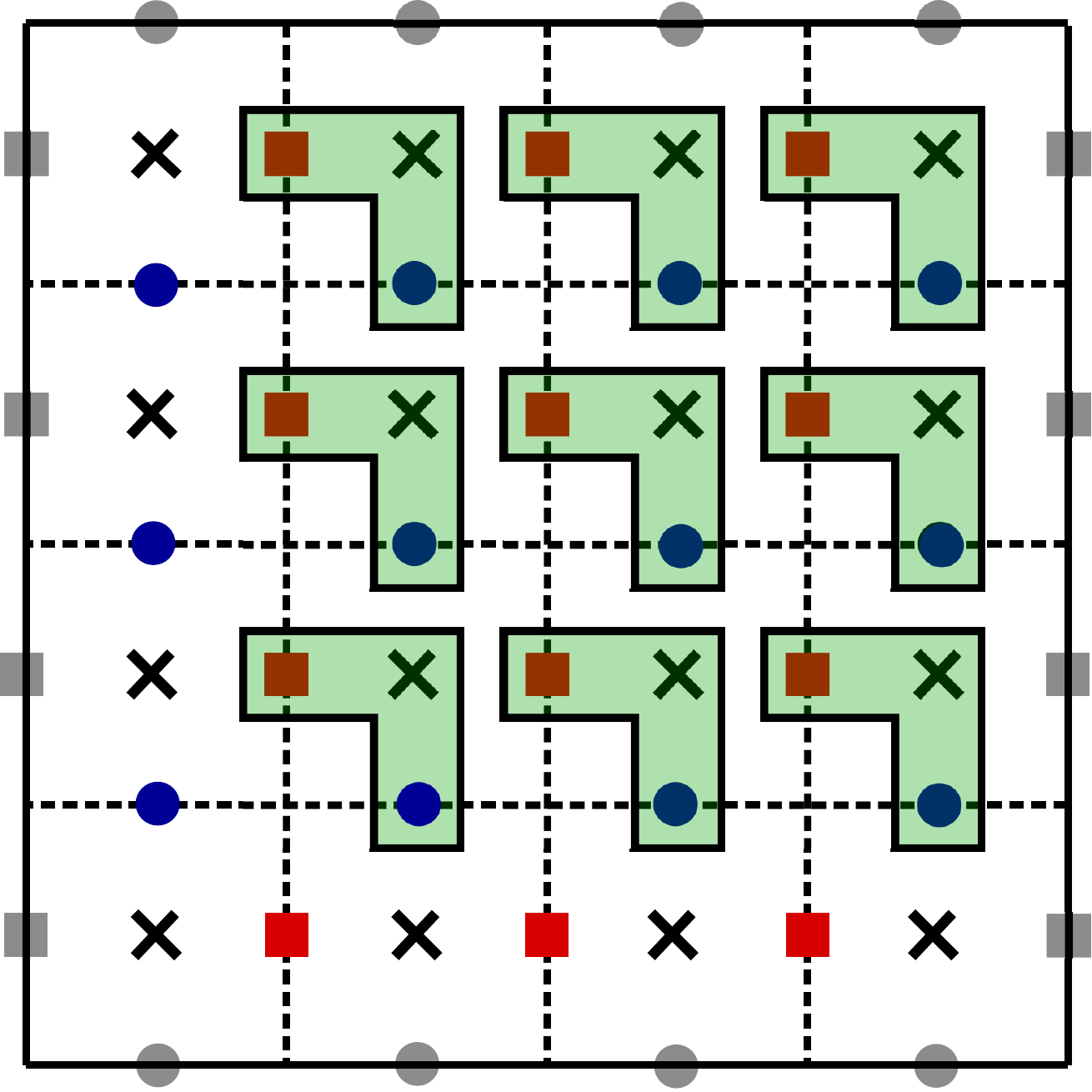}}
\caption{New algorithm for the triad relaxation method. The order of the procedure is performed as visualized in the figures from left to right.}
\label{fig:13}
\end{figure}

The results of the local Fourier analysis, c.f.\ Section \ref{sec:anares}, and the numerical results for periodic boundary conditions, c.f.\ Section \ref{sec4}, in combination with the better parallelelizability and the reduced computational work demonstrate the potential of the triad-wise relaxation method. However, the numerical results for the Stokes system with Dirichlet boundary conditions indicate the issue of the triad-wise method. In what follows, we show how to fix this issue. For this purpose more computational work is needed. The idea is to repeat the relaxation process four times while changing the unknowns contained in one block after each iteration. In addition, we vary the order in which the blocks are updated. This idea is illustated in Figure \ref{fig:13}. 

\begin{figure}[h]
\centering
\includegraphics[width=13 cm]{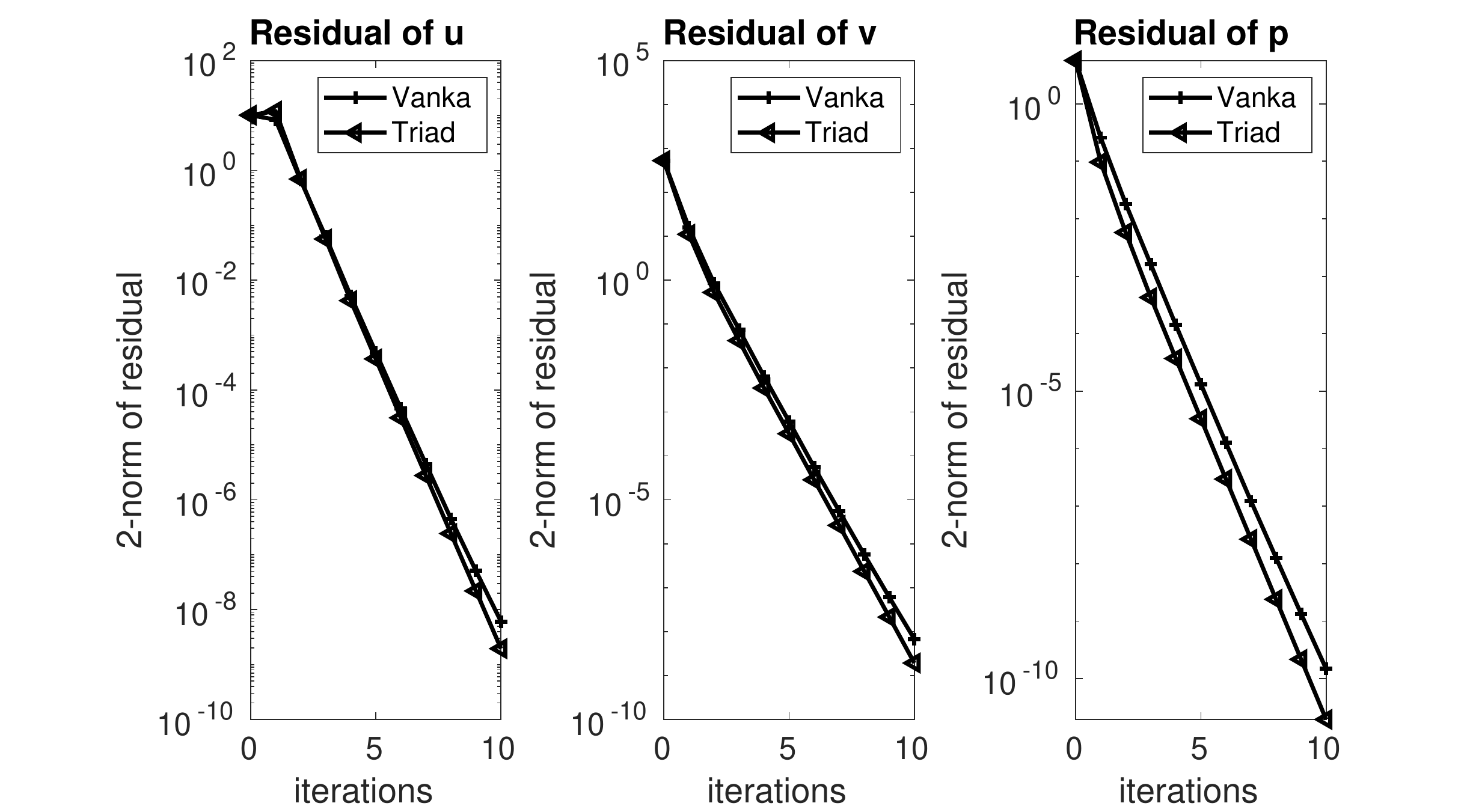}
\caption{Convergence behavior of the multigrid method for the Stokes system with Dirichlet boundary conditions including the new triad algorithm.}
\label{fig:14}
\end{figure}

The order of the iteration is visualized from left to right, i.e., Figure \ref{subfig:13a} illustrates the first iteration and Figure \ref{subfig:13d} the last one. This new algorithm improves the convergence results. 
By changing the order of relaxing different blocks, there is not a significant impact on the convergence factor. However, small variations within a range of 0.03 are noticeable and the relaxation process in lexicographical order, i.e., row-wise starting with the bottom left block, promises best results. The two-grid convergence factor of the updated version of the triad-wise smoother for a homogeneous problem with Dirichlet boundary conditions is $\rho=0.04$. Again, we applied 20 two-grid cycles with two pre- and postsmoothing steps starting with a random initial guess. Results of a multigrid method are visualized in Figure \ref{fig:14}.

The modified procedure causes four times more computational work in comparison with the results of the original version of the triad-wise smoother, c.f. Section \ref{sec4}. This leads to a higher number of arithmetic operations for the triad-wise smoother compared to the Vanka smoother. We have a total number of $139\cdot (N-1)^2$ for the Vanka smoother versus $176\cdot (N-1)^2$ operations for the updated triad smoother. However, the new triad procedure convergences even better than the Vanka smoother. We have $\rho=0.04$ for the triad-wise smoother compared to $\rho=0.10$ for the Vanka smoother. Each of the individual substeps of the modified method can be parallelized in the same way as a whole step of the non-modified method. Nevertheless, the four steps are processed consecutively resulting in more synchronization steps.
%

\section{Parallel implementation aspects}\label{sec6}
So far, we have seen a comparison between an overlapping and a non-overlapping smoother, analyzed their convergence behavior and potential for optimization. Applications for the Stokes equations result in systems of equations that are much larger than the test problems within this paper, i.e.\ with several million unknowns. Consequently, a serial implementation is usually not sufficient. As a result, we compare the parallelization potential of the overlapping smoother with the non-overlapping smoother in the following.

The obvious difference between the two smoothers is the number of unknowns
and equations that are included in one block, while the number of blocks that have to be updated per relaxation is the same for both smoothers. The arithmetic operations that we need for one relaxation step are $139 \cdot N^2$ for the Vanka smoother compared to $44 \cdot N^2$ for the original version of the triad-wise smoother, where $N$ describes the number of blocks. More details regarding the computation of the number of arithmetic operations can be found in \citep{clausphd}. Parallel implementation aspects that we take into account are the required coloring and the communication effort that is needed for computations in parallel.

\begin{figure}[h]
\centering
\subcaptionbox[]{The triad-wise smoother
    \label{subfig:15a}} 
[0.3\textwidth] 
{\includegraphics[width=4 cm]{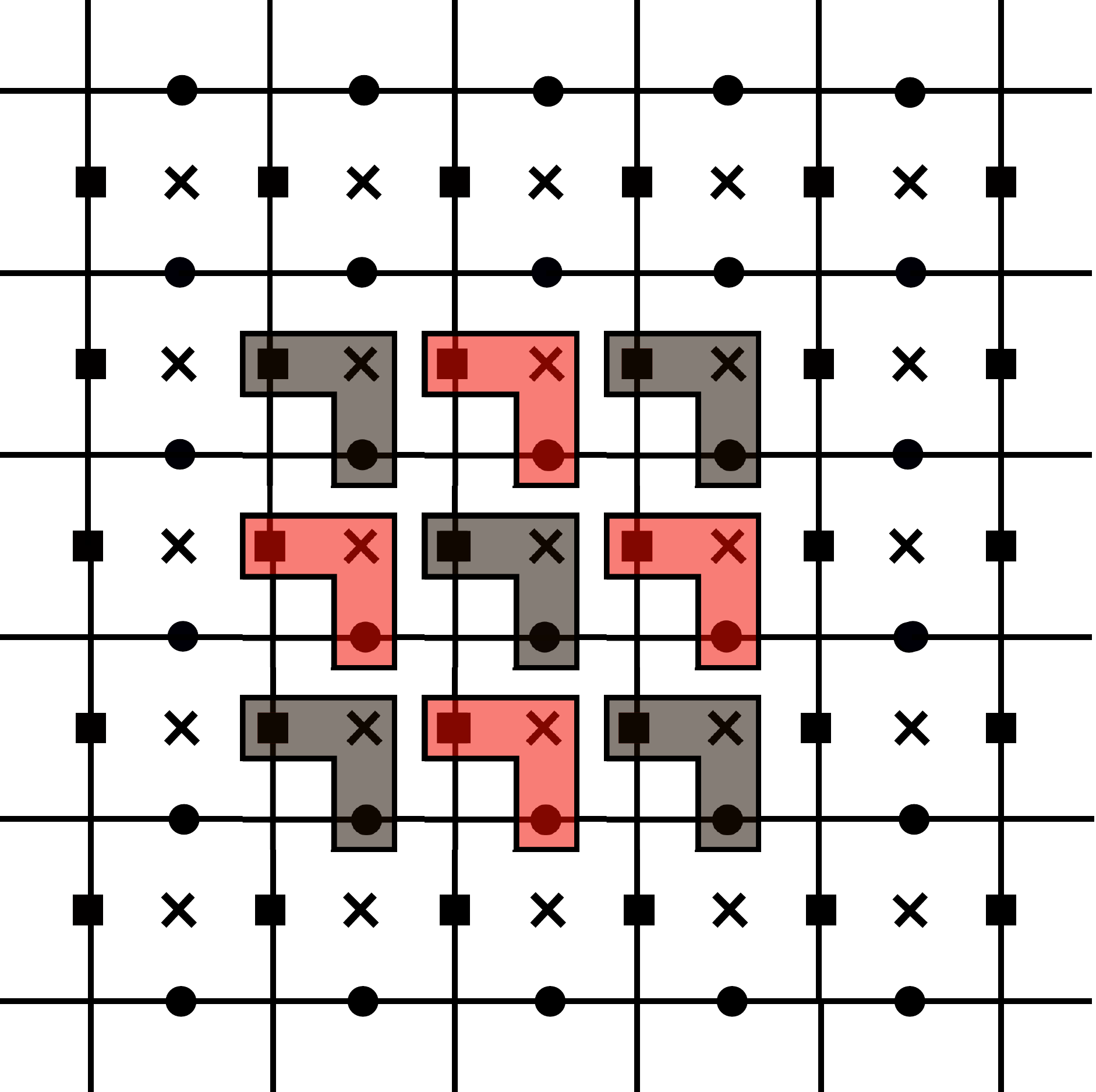}}
\subcaptionbox[]{The Vanka smoother (5 color)
    \label{subfig:15b}} 
[0.3\textwidth] 
{\includegraphics[width=4 cm]{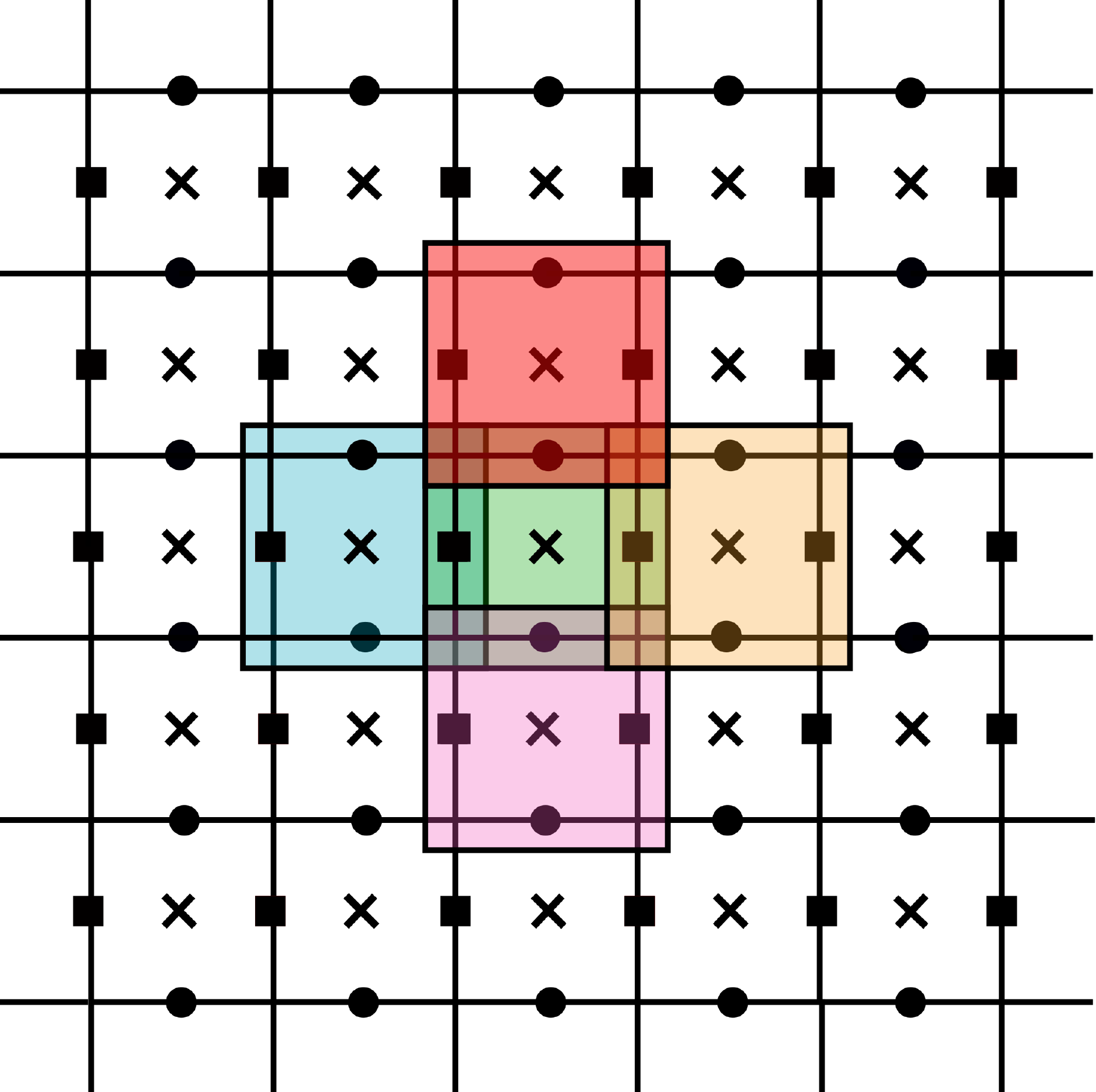}}
\subcaptionbox[]{The Vanka smoother (9 color)
    \label{subfig:15c}} 
[0.3\textwidth] 
{\includegraphics[width=4 cm]{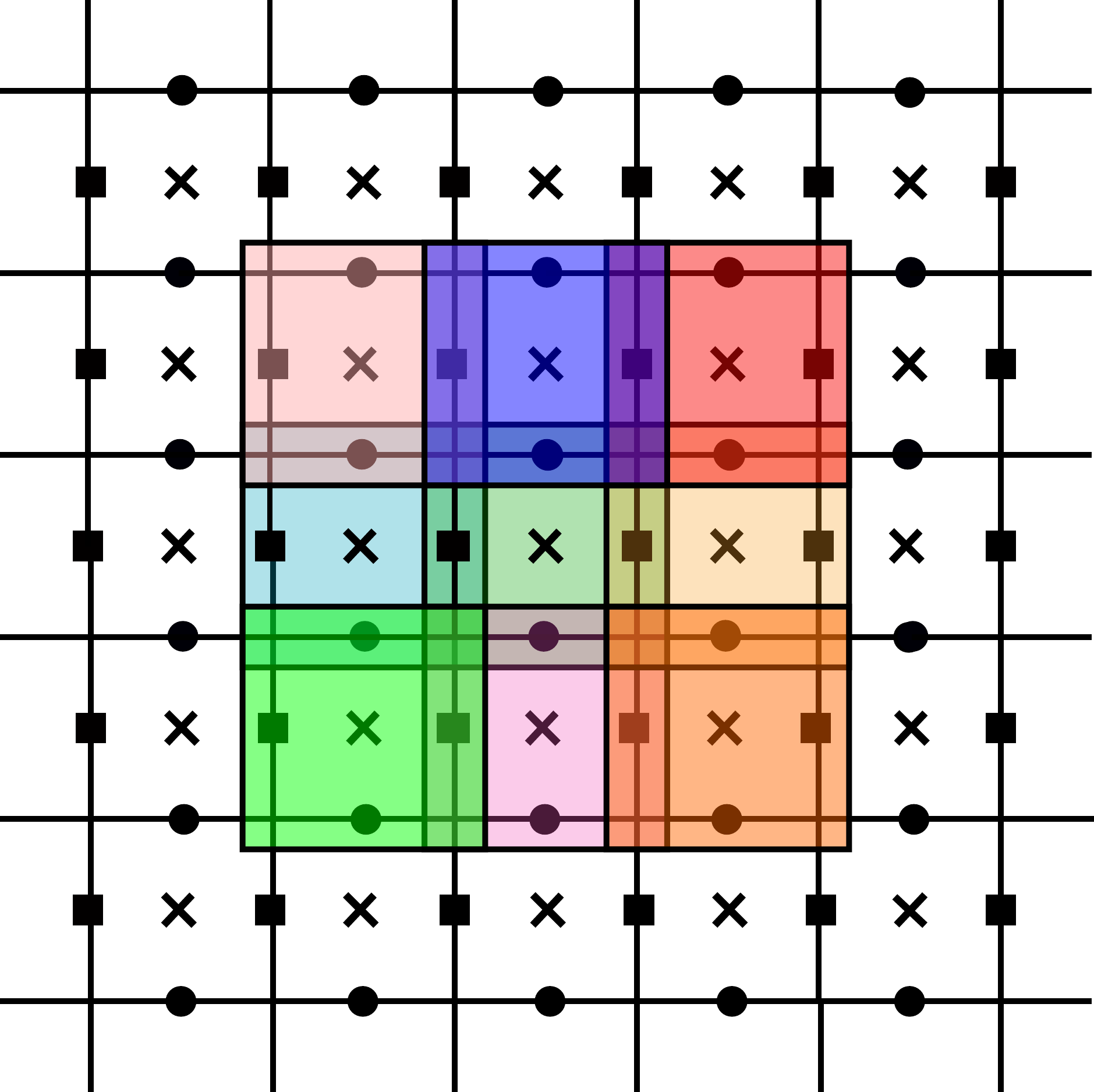}}
\caption[]{Coloring schemes.}
\label{fig:15}
\end{figure}

The required colorings to perform the smoothing steps in parallel are visualized in Figure \ref{fig:15}. The Vanka smoother needs five different colors, while the triad-wise smoother requires a red-black coloring (i.e.\ two colors), only.

\begin{figure}[h]
\centering
\subcaptionbox[]{The 5 color Vanka smoother
    \label{subfig:16a}} 
[0.45\textwidth] 
{\includegraphics[width=6 cm]{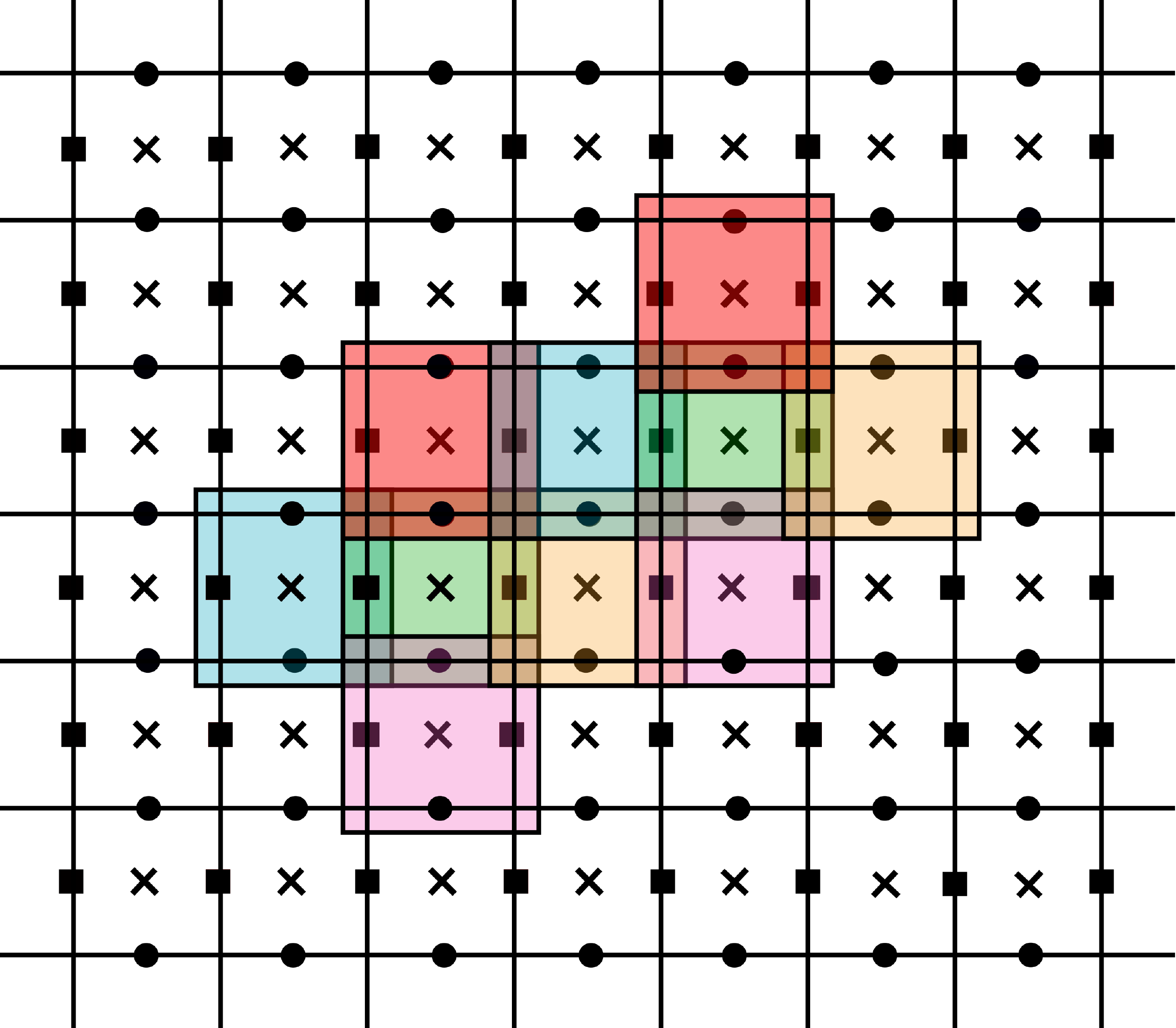}}
\subcaptionbox[]{The 9 color Vanka smoother
    \label{subfig:16b}} 
[0.45\textwidth] 
{\includegraphics[width=6 cm]{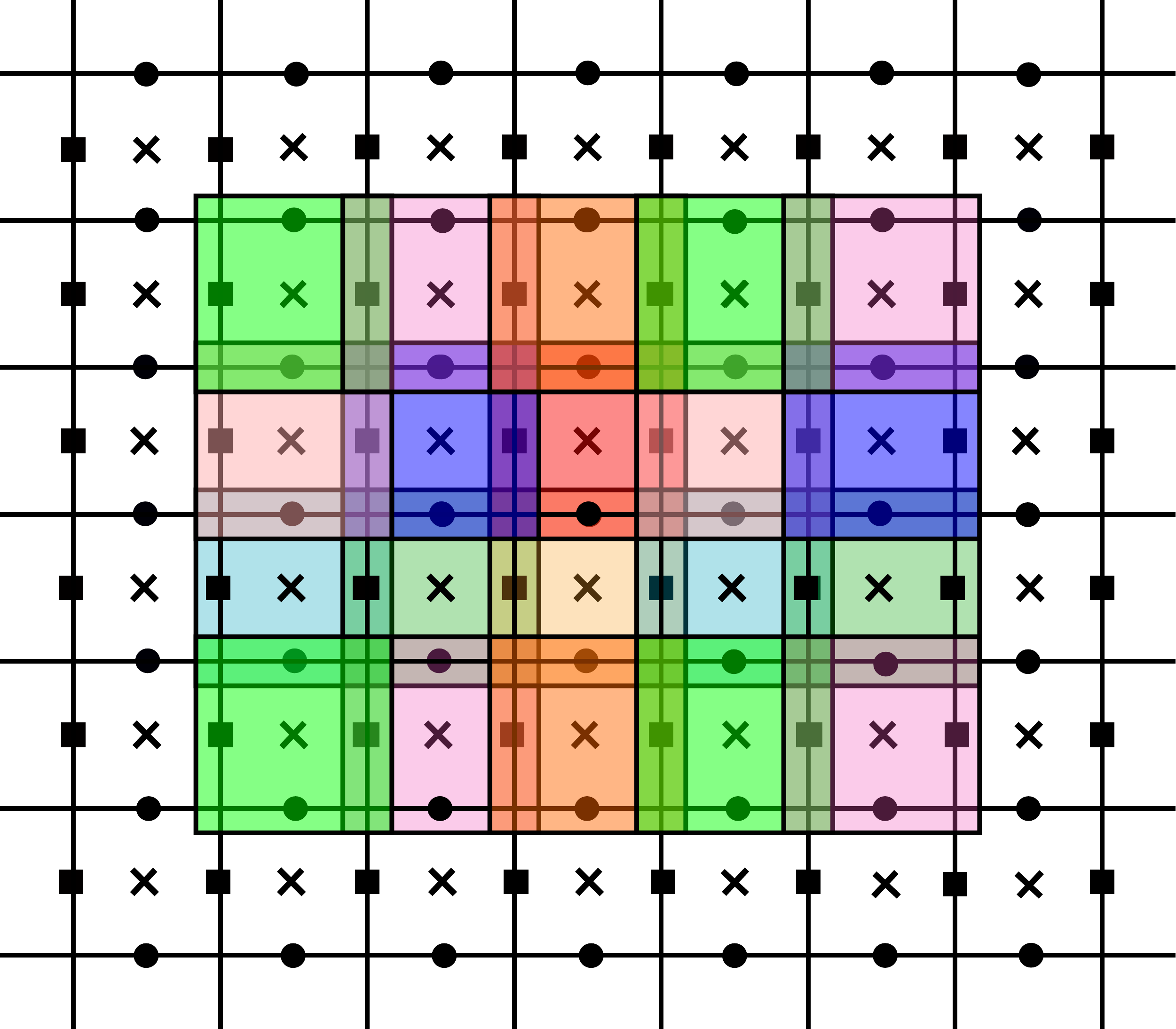}}
\caption{Repetition of the coloring pattern of the Vanka smoother.}
\label{fig:16}
\end{figure}

An important parallelization aspect is the ``dependencies'', i.e., all neighboring unknowns that are needed to update an unknown. To perform a method in parallel without read and write conflicts it is necessary to avoid simultaneous updates of dependent unknowns. The coloring scheme is based on depenencies. The unknowns of one triad block depend on the unknowns of adjacent cells. Consequently, a red-black coloring is sufficient, c.f. Figure \ref{subfig:15b}. That means all unknowns of the red triad blocks can be updated simultaneously while all unknowns of the black triad blocks can be updated simultaneously. Whereas, the coloring scheme of the Vanka smoother is illustrated in Figure \ref{subfig:15a}. Due to the additional unknowns that are included in one Vanka block compared to the triad block, the unknowns do not solely depend on the adjacent blocks. Due to the MAC scheme as given in Equations \eqref{eq:5} - \eqref{eq:7} a five coloring is sufficient to perform the Vanka smoother in parallel.
The repetition of the coloring schemes are visualized in Figure \ref{subfig:15b} for the red-black pattern of the triad-wise smoother and in Figure \ref{fig:16} for the five color Vanka pattern.

The communication effort differs due to the different numbers of unknowns that are updated within each block smoothing step. Moreover, the dependence on updated unknowns of neighboring cells leads to further differences. For the Vanka smoother $60\cdot N^2$ communication steps are necessary, while for the triad-wise smoother $12\cdot N^2$ are needed for one smoothing iteration. One smoothing step of the Vanka smoother requires as many communication steps as five smoothing steps of the triad-wise smoother. These results show the low computational cost in combination with good parallelization properties of the triad-wise smoother. However, a parallel implementation may lead to different convergence factors due to the coloring scheme, i.e. a modified order of updating the unknowns. 

\begin{figure}[h]
\centering
\includegraphics[width=13 cm]{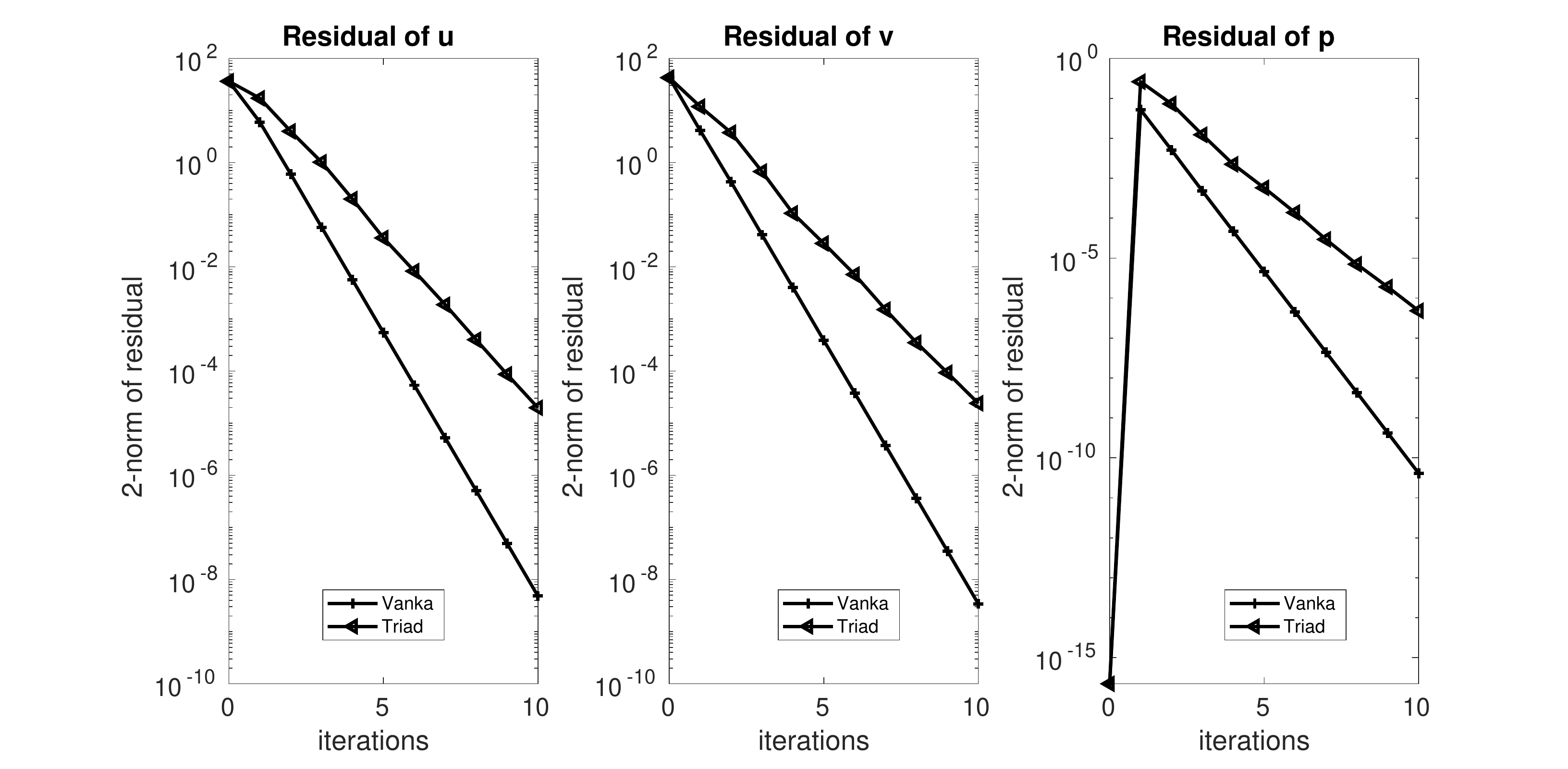}
\caption{Convergence behavior of the multigrid method for the Stokes system with periodic boundary conditions.}
\label{fig:17}
\end{figure}

We implement a two-grid method in serial with the colored versions of both smoothers to get an idea how the coloring influences the convergence factors.
We start with a comparison of numerial results based on the overlapping and non-overlapping smoother for the Stokes equations with periodic boundary conditions. The results are based on the same example and settings as described in Section \ref{sec4}. That leads to convergence factors of $\varphi = 0.29$ for the triad-wise smoother with red-black coloring. Compared to a convergence factor of $\varphi=0.06$ for the 5-color Vanka smoother. We have similar results compared to the serial implementation without a coloring. In consideration of the low computational effort for the triad-wise smoother the convergence factor is reasonable. Multigrid results show the same effect, see Figure \ref{fig:17}.

This is not true for Dirichlet boundary conditions. As we already notice in Section \ref{sec4}, a different ordering of updating the unknowns may change the convergence results especially with Dirichlet boundary conditions. The two-grid method for the colored versions of both smoothers with Dirichlet boundary conditions leads to the following results. The triad-wise smoother with a Gauss-Seidel update of unknowns and a red-black coloring leads to $\varphi=0.58$ as seen in Table \ref{table:tab1}. However, if we implement the modified version of the triad-wise smoother as described in Section \ref{sec5} with a red-black coloring, we get a convergence factor of $\varphi=0.21$.  The novel triad-wise smoother consists of four consecutive steps, in which one single step is similar to the original triad-wise smoother. The consecutive execution of the four steps results in more synchronization steps. In comparison, the 5-color Vanka smoother leads to convergence results of $\varphi=0.99$. This extremely bad convergence rate for the Vanka smoother can be immensely improved by a 9-coloring (c.f. Figure \ref{fig:15} and \ref{fig:16}) because of symmetry that is lost in the 5-color version. We obtain a convergence rate of $\varphi=0.10$ for the two-grid method with the 9-color Vanka smoother. In comparison with the triad-wise smoother, the Vanka smoother with 9 colors clearly leads to better convergence results. As discussed before the computational costs of the Vanka smoother with 5 colors are approximately the same as the costs of the modified triad-wise smoother. The communication costs of the triad-wise smoother increase from $12\cdot N^2$ communication steps to $48\cdot N^2$ due to the repetition of the original triad-wise smoother for four times, which is still below the costs for the Vanka smoother. In addition, the 9-coloring leads to more consecutive processing and therefore decreases the parallelization of the Vanka smoother even more. Using the triad-wise smoother instead of the Vanka smoother compromises the lower convergence rate with slightly lower communication costs and additionally good parallelization properties.

\section{Conclusions}\label{sec7}

Focusing on the results of the serial implementation we draw the following conclusions. For the Stokes equations with periodic boundary conditions the two-grid and multigrid convergence of the well-known Vanka smoother is better than the convergence of the triad-wise smoother.  However, the triad-wise smoother is less expensive than the Vanka smoother.

For the Stokes equations with Dirichlet boundary conditions the good convergence of the Vanka smoother remains unchanged, while the convergence of the triad-wise smoother is not satisfying. However, the modified triad-wise smoother leads to convergence results that are even better than the convergence results of the Vanka smoother. In comparison to the original triad-wise smoother we need to accept less parallelization potential for the modified triad-wise smoother. These results give a good understanding of the challenges arising from different boundary conditions.

For the parallel implementation, we observe similar characteristics. A clear advantage of the triad-wise smoother over the Vanka smoother is given in respect to its low computational costs and its good parallelization properties with a reasonable higher convergence rate. This is especially true for the Stokes equations with periodic boundary conditions. Dirichlet boundary conditions lead to similar difficulties for the parallel implementation compared to the serial algorithm for the triad smoother. In addition, the Vanka smoother needs additional coloring for the two-grid method to obtain a good convergence rate. This leads to a further decrease of parallelization properties of the Vanka smoother and highlights the advantages of our novel triad-wise smoother.

\nocite{*}
\bibliography{wileyNJD-VANCOUVER}%

%

\end{document}